%% file: system.tex
\documentclass[]{interact}

\usepackage{nccmath}
\usepackage{amsfonts}
\usepackage{amsthm}
\usepackage{amsmath}
\usepackage{amscd}
\usepackage{comment}
\usepackage[latin2]{inputenc}
\usepackage{t1enc}
\usepackage{multicol}
\usepackage{multirow}
\usepackage[mathscr]{eucal}
\usepackage{indentfirst}
\usepackage{graphicx}
\usepackage{graphics}
\usepackage{pict2e}
\DeclareMathOperator*{\ess}{ess\,\sup}
\usepackage{epic}
\usepackage{algorithm}
\usepackage{algorithmic}
\usepackage{subfig}
\usepackage{tikz}
\usepackage{float}
\usepackage{caption}
\usepackage{import}
\usepackage{geometry}
\usepackage{epstopdf}% To incorporate .pdf illustrations using PDFLaTeX, etc.
\numberwithin{equation}{section}

\usepackage[numbers,sort&compress]{natbib}% Citation support using natbib.sty
\bibpunct[, ]{[}{]}{,}{n}{,}{,}% Citation support using natbib.sty
\renewcommand\bibfont{\fontsize{10}{12}\selectfont}% Bibliography support using natbib.sty
\newcommand{\vertiii}[1]{{\left\vert\kern-0.25ex\left\vert\kern-0.25ex\left\vert #1 
		\right\vert\kern-0.25ex\right\vert\kern-0.25ex\right\vert}}
	
\theoremstyle{plain}% Theorem-like structures provided by amsthm.sty
\newtheorem{Th}{Theorem}[section]
\newtheorem{Lemma}[Th]{Lemma}

\newtheorem{Prop}[Th]{Proposition}

\theoremstyle{definition}

\theoremstyle{remark}
\newtheorem{Rem}{Remark}

\title{Virtual element method for the system of time dependent nonlinear convection-diffusion-reaction equation}

\author{
\name{M.~Arrutselvi %\textsuperscript{a}
 and E.~Natarajan \thanks{ Corresponding author E.~Natarajan. Email: thanndavam@iist.ac.in}
}
\affil{Department of Mathematics, Indian Institute of Space Science and Technology, Thiruvananthapuram-695547, Kerala, India.}
}
\begin{document}

\maketitle

\begin{abstract} 
In this work, we have discretized a system of time-dependent nonlinear convection-diffusion-reaction equations with the virtual element method over the spatial domain and the Euler method for the temporal interval. For the nonlinear fully-discrete scheme, we prove the existence and uniqueness of the solution with Brouwer's fixed point theorem. To overcome the complexity of solving a  nonlinear discrete system, we define an equivalent linear system of equations. A priori error estimate showing optimal order of convergence with respect to $H^1$ semi-norm was derived. Further, to solve the discrete system of equations, we propose an iteration method and a two-grid method. In the numerical section, the experimental results validate our theoretical estimates and point out the better performance of the two-grid method over the iteration method.
\end{abstract}


\section{Introduction} 
\label{sec:introduction}

Virtual Element Method (VEM) introduced in \cite{vem1} is a generalization of finite element method, that can handle meshes with polygonal elements of both convex and non-convex types. Moreover, implementation of higher order VEM is simple, as the computation of integrals arising in the discrete scheme does not require the explicit knowledge of the basis functions. VEM has been successfully applied to linear and nonlinear partial differential equations, ranging from convection-diffusion problems \cite{vem-our,vem-our2}, parabolic \cite{vem9,vem5}, hyperbolic \cite{vem10,vem6}, elliptic \cite{vem11,vem4,vem7}, and to various problems such as Stokes \cite{vem16}, biharmonic \cite{vem19}, mixed Brinkman \cite{vem21} and linear elasticity \cite{vem22}, $C1$ VEM for Cahn-Hilliard problem \cite{vem15}, 3D elasticity problem \cite{vem14} and so on. 

In this paper, we study VEM for a system of time-dependent nonlinear convection-diffusion-reaction equations that arise in several practical applications. We derive the computable VEM discretization with the help of projection operators in section 3 and then prove the existence and uniqueness of the discrete nonlinear system of equations in section 4. In section 5 we have performed the convergence analysis by deriving error estimates showing optimal rate of convergence. In section 6, we have introduced two-grid method in order to
reduce the computational cost to solve the system of equations. To establish this fact, we have compared the two-grid approach with the standard iterative procedure in the numerical experiments in section 7.
%An outline of the paper is as follows. In section 2, we present the model problem and its  weak formulation, along with  sufficient regularity conditions on the coefficients. In section 3, we introduce the virtual element space and propose a computable nonlinear fully-discrete VEM formulation. The existence and uniqueness of the solution for the discrete scheme is shown in section 4. Since the system are nonlinear in nature, to avoid complications in the approximation, we propose an equivalent linear discrete VEM formulation and  derive an error estimate obtaining optimal order of convergence with respect to $H1$ semi-norm, in section 5. In section 6, we propose two methods, namely, iterative method and two-grid method, to solve the arising system of equations. Results of the numerical experiments conducted on two examples, are presented in section 7.

\subsection{Notations}
Let $D$ be a subset of $\mathbb{R}^2$ with boundary $\partial D$. The space $L^2(D)$ consists of square integrable functions with inner product $(u,v)_{0}:=\int_Duv\,dx$ and norm $\|\cdot\|_{0}:=(\cdot,\cdot)^{\frac{1}{2}}_{0}$.  
For index $\,s\in\mathbb{N},\,$ $H^s(D)$ is the usual Sobolev space with norm $\|u\|_{s}^2:=\sum\limits_{\alpha\leq s}\|\frac{\partial^\alpha u}{\partial x^\alpha}\|^2_{0}$. 
%The space $H^1_0(D)\subset H^1(D)$ contains functions with zero value on $\partial D$ and the space $H^{-1}(D)$ is the dual of $H^1_0(D)$.
Let $X$ be a space with inner product $(\cdot,\cdot)_X$ and norm $\|\cdot\|_X$. We define $L^p(0,T;X),\,\,1\leq p\leq \infty$,\, be space of measurable functions  $u:[0,T]\rightarrow X$ with norms $\|u\|_{L^p(0,T;X)}$ defined as
\begin{eqnarray}
\|u\|_{L^p(0,T;X)}^p:=\textstyle \int_0^T\,\|u(t)\|_X^p\,\quad 1\leq p<\infty,\qquad
\|u\|_{L^\infty(0,T;X)}:=  \ess_{0\leq t\leq T} \|u(t)\|_X.\nonumber
\end{eqnarray}
%Similarly $C(0,T;X)$ denote the set of continuous functions  $u:[0,T]\rightarrow X$.
For $k\in \mathbb{N}$,\,\,we denote \, $[X]^k:=\,X\times X\times...\times X\,\,(k- \text{times})$. The space $[X]^k$ is endowed with the inner product $\boldsymbol{\langle u,v\rangle}:=\sum\limits_{i=1}^k (u_i,v_i)_X$,\, and the induced norm $\boldsymbol{\|u\|}:=\boldsymbol{\langle u,v\rangle}^\frac{1}{2}$.
%\, Then any\, $\boldsymbol{u}\in \bigoplus\limits_{i=1}^kX$ is of the form $\boldsymbol{u}=(u_1,...,u_k) $ with $u_i\in X,\,\,i=1,...,k.$ The space $\bigoplus\limits_{i=1}^kX$ is endowed with the inner product $\boldsymbol{\langle u,v\rangle}:=\sum\limits_{i=1}^k (u_i,v_i)_X$,\, and the induced norm $\boldsymbol{\|u\|}:=\boldsymbol{\langle u,v\rangle}^\frac{1}{2}$.

\section{Model Problem}
\label{sec:modelproblem}
Let $\Omega\subset \mathbb{R}^2$ be a bounded domain with piecewise smooth boundary, time interval $(0,T)\subset \mathbb{R}$, $m\in \mathbb{N}$ and $\mathsf{D}:=(0,T)\times\Omega$. We consider the following system of $m$ time-dependent convection-diffusion-reaction equations :
\begin{eqnarray}
&&\textstyle\dfrac{\partial u_i}{\partial t}-\xi_i\Delta u_i+(\omega_1,\omega_2)\cdot\nabla u_i+u_i\sum\limits_{j=1}^m\mathsf{A}(i,j)u_j+\sum\limits_{\ell,j\neq i}\mathsf{Q}(i,\ell,j)u_\ell u_j\nonumber\\
&&\qquad+\textstyle\sum\limits_{j=1}^m\mathsf{R}(i,j)u_j\, =\, f_i(t,x)\qquad \text{in} \,\,\, \mathsf{D},\,\quad i=1,2,...,m.\quad\label{modelprob}\\
&& u_i(t,x)\,= \,0\,\,\,\qquad\qquad \,\,\text{on}\,\,\,(0,T)\times\partial\Omega,\quad i=1,...,m.\quad\label{bc}\\
&& u_i(0,x)\, =\, u_i^0(x) \qquad\quad \text{in}\,\,\,\Omega,\,\quad i=1,...,m.\quad\label{ic}
\end{eqnarray}
Here, the convection parameter $(\omega_1,\omega_2)\in [L^2(\,0,T;(H^1(\Omega)\cap L^\infty(\Omega))\,)]^2$ along  with $(\nabla\cdot(\omega_1,\omega_2))(t,x)=0$ a.e in $\mathsf{D}$ and for $i=1,...,m\,$ we have
$\xi_i\in L^2(\,0,T; L^\infty(\Omega)\,)$ are bounded positive diffusion parameters,\,the source or sink functions $f_i\in L^2(\,0,T; L^2(\Omega)\,)$. The kinetic coefficients $\mathsf{A}$,\,\,$\mathsf{Q}$,\,\,$\mathsf{R}$ are assumed to be sufficiently regular functions of $x$ and $t$, to ensure an unique solution for \eqref{modelprob}-\eqref{ic}.  
The homogeneous Dirichlet boundary condition \eqref{bc} was assumed for simplifying the analysis presentation. Our numerical formulation and analysis with slight modification can handle more general boundary conditions arising in practical applications involving combinations of Dirichlet and Neumann boundary conditions.    

Let us denote $\boldsymbol{\mathcal{H}}:=\big[\,L^2(0,T;H_0^1(\Omega))\cap C(0,T;L^\infty(\Omega))\,\big]^m.$ Then the variational formulation for the system \eqref{modelprob}-\eqref{ic} is : \,
Find $\boldsymbol{u}\in \boldsymbol{\mathcal{H}}$ with ${\partial u_i}/{\partial t}\in L^2(0,T;H^{-1}(\Omega))\,\,\,(i=1,...,m)\,$ such that
\begin{eqnarray}
&& u_i(0,x)\, =\, u_i^0(x) \qquad\quad \text{in}\,\,\,\Omega,\,\quad i=1,...,m,\nonumber\\
&&\textstyle \sum\limits_{i=1}^m\Big\lbrace\,\left(\frac{\partial u_i}{\partial t}, v_i\right)_{0}+\left(\, \xi_i\nabla u_i,\nabla v_i\,\right)_{0}+\left( (\omega_1,\omega_2)\cdot\nabla u_i,\,v_i\, \right)_{0}+\Big( u_i\big(\textstyle\sum\limits_{j=1}^m\mathsf{A}(i,j)u_j\big),\,v_i\, \Big)_{0}\nonumber\\
&&\quad +\Big( \big(\textstyle\sum\limits_{\ell,j\neq i}\mathsf{Q}(i,\ell,j)u_\ell u_j\big),\,v_i\, \Big)_{0}+\Big( \big(\textstyle\sum\limits_{j=1}^m\mathsf{R}(i,j)u_j\big),\,v_i\, \Big)_{0}\,\Big\rbrace = \textstyle \sum\limits_{i=1}^m (\,f_i,\,v_i\,)_0, \label{weakform}\\
&&\qquad\qquad\qquad\qquad\qquad\qquad\qquad\qquad\qquad\qquad \forall \,t\in [0,T],\quad \forall\,\boldsymbol{v}\in [H_0^1(\Omega)]^m.\nonumber
\end{eqnarray}

We assume that the exact solution $\boldsymbol{u}$ of \eqref{modelprob}-\eqref{ic} satisfies the following :
\begin{eqnarray}
    \exists \,\mu>0\quad \text{such that}\quad \|u_i\|_\infty\leq \mu,\quad\forall i\in\lbrace 1,2,...,m\rbrace,\quad\forall \,t\in [0,T].\label{assump}
\end{eqnarray}

\section{The Virtual Element Method}
\label{sec:vemspace}
In this section we brief the definition of the virtual element space. Let $\mathbb{P}_p(E)$ denote the space of polynomials with degree  less than or equal to $p$ on $E$.

We consider $\lbrace \mathcal{T}_h \rbrace_{h>0} $ to be  a family of  partitioning of $\Omega$ consisting of polygonal elements. For some theoretical estimates to hold, we make the following assumptions on every element in $\mathcal{T}_h$.
For an element $E\in \mathcal{T}_h$, let $h_E$ be the element diameter and
 $h$ represent the maximum diameter over all $E\in\mathcal{T}_h$.
In addition, for constants $\delta>0,\,c>0\,$  independent of $h$,\,$E,$\, we assume an element $E$ of $\mathcal{T}_h$ fulfil the following (see \cite{vem1}) : 
(i)  $E$ is star-shaped with respect to a disc $D_{\delta}$ of radius $\delta \, h_E$, \,
(ii) for any edge $e \subset E$, the length $|e| \geq c \, h_E$,\,and\,(iii) boundary of $E$ is made up of finite number of edges.

 Let us define some necessary operators that projects elements in $H^1(E)$ onto $\mathbb{P}_p(E)$. We define the the $L^2$ projection operator  $\mathcal{P}_p^{0}$ (see \cite{vem4}) by,
\begin{equation}
\left( u-\mathcal{P}_p^0 u,\upsilon_p \right)_{E} = 0 \quad \forall\, \upsilon_p \in \mathbb{P}_p(E).
\label{proj2}
\end{equation}

The gradient projection operator $\mathcal{P}_p^{\nabla}: H^1(E) \rightarrow \mathbb{P}_p(E)$ is such that (see \cite{vem4}),
\begin{equation}
\left(\nabla (u-\mathcal{P}_p^{\nabla}u), \nabla \upsilon_p  \right)_{E} = 0 \quad \forall \,\upsilon_p \in \mathbb{P}_p(E)\, \,\,\text{and} \,\,
{\int_{\partial E} (\mathcal{P}_p^{\nabla}u-u) \, ds = 0.}
%P^0(\Pi_p^{\nabla}u-u) = 0
\label{proj1}
\end{equation}

For $u\in H^1(E)$, we shall determine  $\boldsymbol{\mathcal{P}}_{p-1}^0(\nabla u) \in (\mathbb{P}_{p-1}(E))^2$ by,
	\begin{equation}
	\left( \nabla u-\boldsymbol{\mathcal{P}}_{p-1}^0 \nabla u,\boldsymbol{\upsilon}_{p-1} \right)_{E} = 0 \quad \forall\, \boldsymbol{\upsilon}_{p-1} \in (\mathbb{P}_{p-1}(E))^2.
	\label{proj3}
	\end{equation}

Consider the auxiliary space $W_E^p$ (see \cite{vem52}) for each $E \in \mathcal{T}_h$ by,
\begin{equation*}
W_E^p=\left \lbrace v \in H^1(E) \cap C^0(\partial E): v_{|e} \in \mathbb{P}_p(e) \, \forall \, \text{edge} \, { e \in \partial E,}\, \Delta v \in \mathbb{P}_p(E) \right \rbrace.
\end{equation*}

\noindent Let $\left(\mathbb{P}_p(E)/\mathbb{P}_{p-2}(E)\right)$ be the set of polynomials of degree exactly equal to $p-1$ and $p$.
Now we define the local virtual element space $\mathcal{V}_E^p$ as follows,
\begin{equation}
\mathcal{V}_E^p=\left \lbrace\, u \in W_E^p \quad \mathrm{s.t.} \, \, \left(u-\mathcal{P}_p^{\nabla}u,q \right)_{E}=0  \quad \forall q \in \left(\mathbb{P}_p(E)/\mathbb{P}_{p-2}(E)\right)\, \right \rbrace. \end{equation}
We consider the 
following set of degrees of freedom  on $V_E^p$ ,\, ($V_1$) the values of $u$ at the $n(E)$ vertices of polygon $E$,\, ($V_2$) the values of $u$ at $(p-1)$ internal Gauss-Lobatto quadrature nodes of every edge $e \in \partial E$,
\, and ($V_3$) the moments up to order $p-2$ of $u$ in $E$, i.e.,
\begin{equation*}
V_3\,:=\,\big\lbrace\,(\, u, \, \upsilon_{p-2} \,)_E\,\,:\,\, \upsilon_{p-2} \in \mathbb{P}_{p-2}(E)\,\big\rbrace.
\end{equation*}

%\begin{figure}[H]
%	\centering 
%	\includegraphics[width=13cm]{vemsc-pics/dof_pic.pdf}	
%	\caption{Degrees of freedom for $p=1,2,3$ (from left to right). We represent $V_1$ by circle, $V_2$ by
%		triangle and the moments $V_3$ by solid rectangle.}
%	\label{degfm}
%\end{figure}

We note that $(V_1)-(V_3)$ determine  $u\in \mathcal{V}_E^p$ uniquely on the polygon $E$, (see \cite{vem52}) . Now we
define the global virtual element space $V_h^p$ by,
\begin{equation}
V_h^p = \lbrace  \,\,\,u \in H_0^1(\Omega) \,\,\, :\, \, u_{|E} \in V_E^p \, \, \,\,\, \forall E \,\in \mathcal{T}_h  \,\,\,\rbrace. \label{wj}
\end{equation}

\subsection{VEM formulation}
Consider a partition of $(0,T)$ into $N$ disjoint intervals of width $\delta t\,(=T/N)$ and denote $t_n=n\,\delta t,\,\,n=0,1,...,N.$ \,Let us denote by $U_i^n$, the function $U_i$ associated to time $t_n$, i.e. $U_i^n:=U_i(t_n)\in V_h^p$.
We use backward Euler method to discretize the time derivative term. Hereafter we use the inner product $(\cdot,\cdot)_0$ and norms $\|\cdot\|_k$ split over the mesh elements in $\mathcal{T}_h$, i.e, $(\cdot,\cdot)_0:=\sum\limits_{E\in\mathcal{T}_h}(\cdot,\cdot)_{0,E},\,\,$ and $ \,\|\cdot\|_k:=\sum\limits_{E\in\mathcal{T}_h}\|\cdot\|_{0,E}$.% and  denote $\boldsymbol{\mathcal{V}_h^p}:=\bigoplus\limits_{i=1}^m V_h^p.$ 

As the functions in $V_h^p$ are not explicitly known, using the projection operators described in section \ref{sec:vemspace}, we define a $computable$ fully discrete virtual element scheme equivalent to \eqref{weakform} as : 

\noindent for each $n=1,2,...,N,$\, Find $(U_1^n,U_2^n,...,U_m^n)\in \boldsymbol{\mathcal{V}_h^p}$ such that
\begin{eqnarray}
&& U^0_i(x)\, =\, u_i^0(x) \qquad\quad \text{in}\,\,\,\Omega,\,\quad i=1,...,m,\nonumber\\
&&\textstyle \sum\limits_{i=1}^m\Big\lbrace\,\left(\mathcal{P}(\, U_i^n-U_i^{n-1}\,),\, \mathcal{P}V_i\right)_{0}+S_1((I-\mathcal{P})\,(U_i^n-U_i^{n-1}),\,(I-\mathcal{P})V_i)\nonumber\\
&&\quad+(\delta t)\Big(\, \xi_i\boldsymbol{\mathcal{P}}\nabla U^n_i,\boldsymbol{\mathcal{P}}\nabla V_i\,\Big)_{0}+(\delta t) S_2((I-\Pi_p^\nabla)U_i^n,(I-\Pi_p^\nabla)V_i)+(\delta t)\Big( (\omega_1,\omega_2)\cdot\boldsymbol{\mathcal{P}}\nabla U^n_i,\,\mathcal{P}V_i\, \Big)_{0}\nonumber\\
&&\quad+(\delta t)\Big( \mathcal{P}U^n_i\big(\textstyle\sum\limits_{j=1}^m\mathsf{A}(i,j)\mathcal{P}U^n_j\big),\,\mathcal{P}V_i\, \Big)_{0} +(\delta t)\Big( \big(\textstyle\sum\limits_{\ell,j\neq i}\mathsf{Q}(i,\ell,j)\mathcal{P}U^n_\ell \mathcal{P}U^n_j\big),\,\mathcal{P}V_i\, \Big)_{0}\nonumber\\
&&\quad+(\delta t)\Big( \big(\textstyle\sum\limits_{j=1}^m\mathsf{R}(i,j)\mathcal{P}U^n_j\big),\,\mathcal{P}V_i\, \Big)_{0}\,\Big\rbrace =  (\delta t)\sum\limits_{i=1}^m (\,f^n_i,\,\mathcal{P}V_i\,)_0\qquad \forall\,\boldsymbol{V}\in \boldsymbol{\mathcal{V}_h^p}, \label{vemform}
\end{eqnarray}
 where  $S_1$\,and $S_2$ are symmetric bilinear functions detailed below.
 
 Consider the local symmetric bilinear maps $\widetilde{S_1^E}$ and $\widetilde{S_2^E}$ defined on $V_h^E\times V_h^E$ such that there exists non-zero positive constants $\alpha_{\ast}$, $\alpha^{\ast}$, $\beta_{\ast}$ and $\beta^{\ast}$, with $\alpha_{\ast}\leq \alpha^{\ast}$ and $\beta_{\ast}\leq \beta^{\ast}$, independent of $h_E$, and $\forall u_h, v_h\in V_E^p\setminus \mathbb{P}_p(E),$
\begin{eqnarray}
\alpha_{\ast}( u_h,\, v_h)_{0,E}&\leq& \widetilde{S^E_1}(u_h,v_h) \leq \alpha^{\ast}(u_h,\,v_h)_{0,E},
\label{vemstab-se1}\\
\beta_{\ast}(\nabla u_h, \nabla v_h)_{0,E}&\leq& \widetilde{S^E_2}(u_h,v_h) \leq \beta^{\ast}(\nabla u_h, \nabla v_h)_{0,E} .
\label{vemstab-se2}
\end{eqnarray}
A variety of computable choices for $ \widetilde{S^E_1}$ and $ \widetilde{S^E_2}$ can be found in the literature. Let us choose $S_1$ and $S_2$, as follows,

\begin{eqnarray}
\textstyle
S_1(u_h,v_h)=\sum\limits_{E\in\mathcal{T}_h}\widetilde{S^E_1}(u_h,v_h)\qquad\text{and}\qquad S_2(u_h,v_h)=\sum\limits_{E\in\mathcal{T}_h}(\Bar{\xi_i}^E)\,\widetilde{S^E_2}(u_h,v_h),\qquad\label{stabliomega}
\end{eqnarray}
where $\Bar{\xi_i}^E:=\sup\limits_{x\in E}|\xi_i(x,t_n)|$. In the sequel, we adopt the following notations,
\begin{eqnarray}
m_h(u_h,v_h)&:=& \left(\mathcal{P} u_h,\, \mathcal{P}v_h\right)_{0}+S_1((I-\mathcal{P})u_h,(I-\mathcal{P})v_h),\quad\label{massmat}\\
a_h(u_h,v_h)&:=& \big(\, \xi_i\boldsymbol{\mathcal{P}}\nabla u_h,\boldsymbol{\mathcal{P}}\nabla v_h\,\big)_{0}+ S_2((I-\Pi_p^\nabla)u_h,(I-\Pi_p^\nabla)v_h).\label{stiffmat}
\end{eqnarray}

\section{Theoretical Analysis}
Here, we prove the existence and uniqueness of solution for the fully discrete scheme \eqref{vemform}. First, we state a variant of Brouwer's fixed point that will be used to prove the existence of a solution for \eqref{vemform}.

\begin{Prop}
Let H be a finite dimensional Hilbert space with inner product $\langle \cdot, \cdot \rangle_H$ and norm $\|\cdot\|_H$. Let $L:H\rightarrow H$ be a continuous map. 
If there exists \, $\rho>0$ such that, 
$\langle L(w), w\rangle_H > 0, \, \forall w\in H$ with $\|w\|_H = \rho,$\,then there exists a $z\in H$ such that $L(z)=0$ and $\|z\|_H\leq \rho.$
\label{prop-Brouwer}
\end{Prop}

In our analysis, we assume (later proved in sec. \ref{sec:linearscheme}) that any solution $\boldsymbol{U^{n-1}}$ of \eqref{vemform} satisfies the following :
\begin{eqnarray}
\exists\,\mu>0\quad\text{s.t}\quad \|U^{n-1}_i\|_\infty\leq (\mu+1)\qquad \forall\, i\in \lbrace 1,...,m\rbrace.\label{assump1}
\end{eqnarray}

\begin{Rem}
	On a bounded Lipschitz element $E$, for $\Pi_k^0 w\in \mathbb{P}_k(E)\subset H^1(E)$,  using the compact Sobolev embedding $H^1(E)\hookrightarrow L^p(E),\,\,2\leq p<\infty$ and (2.44) in \cite{vemsuszane} , we have
	\begin{eqnarray}
		\|\Pi_k^0w\|_{L^p(E)}&\leq&  C\,\|\Pi_k^0w\|_{1,E}\leq  C\,\|w\|_{1,E}. \label{sobolevimbedding} 
	\end{eqnarray}
\end{Rem}

\begin{Th}
Under the assumption \eqref{assump1} and 
for sufficiently small $(\delta t)$, we assume the existence of unique tuple  $\boldsymbol{U^k}:=(U_1^{k},U_2^k,...,U_m^k)\in\boldsymbol{\mathcal{V}_h^p}$ for $0\leq k\leq n-1$ satisfying  \eqref{vemform}. Then for all $1\leq n\leq N$ there exists a unique solution $\boldsymbol{U^n}:=(U_1^n,U_2^n,...,U_m^n)\in\boldsymbol{\mathcal{V}_h^p}$ to the discrete problem \eqref{vemform}.
\label{thm:nonlinear}
\end{Th}
\begin{proof}
We observe that $\boldsymbol{\mathcal{V}_h^p}$ is a finite dimensional Hilbert space with inner product $\boldsymbol{\langle W,V\rangle}=\sum\limits_{i=1}^m (\nabla W_i,\nabla V_i)_0$ and the induced norm $\boldsymbol{\|W\|}=\Big(\,\sum\limits_{i=1}^m\|\nabla W_i\|_0^2\,\Big)^\frac{1}{2}$.

For $1\leq n\leq N$, let us define a map $\boldsymbol{L}:\boldsymbol{\mathcal{V}_h^p}\rightarrow\boldsymbol{\mathcal{V}_h^p}$ such that
\begin{eqnarray}
&&\boldsymbol{\langle L(W),\,V\rangle}=\textstyle
\sum\limits_{i=1}^m\bigg\lbrace\, m_h(W_i^n,V_i)-m_h(U_i^{n-1},V_i) +(\delta t)\, a_h(W_i^n,V_i)\nonumber\\
&&\quad+(\delta t)\Big( (\omega_1,\omega_2)\cdot\boldsymbol{\mathcal{P}}\nabla W^n_i,\,\mathcal{P}V_i\, \Big)_{0}+(\delta t)\Big( \mathcal{P}W^n_i\big(\textstyle\sum\limits_{j=1}^m\mathsf{A}(i,j)\mathcal{P}W^n_j\big),\,\mathcal{P}V_i\, \Big)_{0} \nonumber\\
&&\quad+(\delta t)\Big( \big(\textstyle\sum\limits_{\ell,j\neq i}\mathsf{Q}(i,\ell,j)\mathcal{P}W^n_\ell \mathcal{P}W^n_j\big),\,\mathcal{P}V_i\, \Big)_{0}+(\delta t)\Big( \big(\textstyle\sum\limits_{j=1}^m\mathsf{R}(i,j)\mathcal{P}W^n_j\big),\,\mathcal{P}V_i\, \Big)_{0}-(\delta t) (\,f^n_i,\,\mathcal{P}V_i\,)_0\,\bigg\rbrace.\nonumber\\
&&\qquad\qquad\qquad \equiv \textstyle\sum\limits_{i=1}^m \mathcal{F}_i(\boldsymbol{W},\boldsymbol{V}).\qquad\label{eq21}
\end{eqnarray}
We can prove that $\boldsymbol{L}$ is well defined with the help of Riesz representation theorem. Now we shall show that $\boldsymbol{L}$ is a continuous mapping. Let $\epsilon>0$ be given. Consider arbitrary $\boldsymbol{W}\in \boldsymbol{\mathcal{V}_h^p}$, choose a $\nu>0$ (depending on $\epsilon,\,\boldsymbol{W}$; for existence of such $\nu$, see Remark \ref{rem:continuity}) and let $\boldsymbol{Z}\in\boldsymbol{\mathcal{V}_h^p}$ be any tuple such that $\boldsymbol{\|W-Z\|}<\nu$.  Denote $\boldsymbol{\psi}=\boldsymbol{W-Z} $ and $\boldsymbol{L^\psi}:=\boldsymbol{L(W)-L(Z)}$. Then
\begin{eqnarray}
\boldsymbol{\|L(W)-L(Z)\|}^2=\boldsymbol{\langle L^\psi,\,L^\psi\rangle}&=&\boldsymbol{\langle L(W),\,L^\psi\rangle}-\boldsymbol{\langle L(Z),\,L^\psi\rangle}\nonumber\\
&\equiv&\textstyle \sum\limits_{i=1}^m \big\lbrace\, \mathcal{F}_i(\boldsymbol{W},\boldsymbol{L^\psi})-\mathcal{F}_i(\boldsymbol{Z},\boldsymbol{L^\psi})\,\big\rbrace.\label{eq3}
\end{eqnarray}

Now we estimate the terms in \, $I:=\mathcal{F}_i(\boldsymbol{W},\boldsymbol{L^\psi})-\mathcal{F}_i(\boldsymbol{Z},\boldsymbol{L^\psi})$.

\begin{eqnarray}
I&=&m_h(\psi_i^n,L^\psi_i)+(\delta t)\, a_h(\psi_i^n,L^\psi_i)+(\delta t)\Big( (\omega_1,\omega_2)\cdot\boldsymbol{\mathcal{P}}\nabla \psi^n_i,\,\mathcal{P}L^\psi_i\, \Big)_{0}\nonumber\\
&&\quad+(\delta t)\Big( \mathcal{P}W^n_i\big(\textstyle\sum\limits_{j=1}^m\mathsf{A}(i,j)\mathcal{P}W^n_j\big)-\mathcal{P}Z^n_i\big(\textstyle\sum\limits_{j=1}^m\mathsf{A}(i,j)\mathcal{P}Z^n_j\big),\,\mathcal{P}L^\psi_i\, \Big)_{0} \nonumber\\
&&\quad+(\delta t)\Big( \big(\textstyle\sum\limits_{\ell,j\neq i}\mathsf{Q}(i,\ell,j)\big[\,\mathcal{P}W^n_\ell \mathcal{P}W^n_j-\mathcal{P}Z^n_\ell \mathcal{P}Z^n_j\,\big]\,\big),\,\mathcal{P}L^\psi_i\, \Big)_{0}\nonumber\\
&&\quad+(\delta t)\Big( \big(\textstyle\sum\limits_{j=1}^m\mathsf{R}(i,j)\,\mathcal{P}\psi^n_j\,\big),\,\mathcal{P}L^\psi_i\, \Big)_{0}:=\sum\limits_{r=1}^6I_r.\quad\label{eq4}
\end{eqnarray}
Using the definition \eqref{massmat}, property \eqref{vemstab-se1} and Cauchy-Schwarz inequality, we get
\begin{eqnarray}
I_1:=m_h(\psi_i^n,L^\psi_i)\leq \max\lbrace1,\alpha^*\rbrace\,(\psi_i^n,\,L^\psi_i)_0 \leq \max\lbrace1,\alpha^*\rbrace\,\|\psi_i^n\|_0\,\|L^\psi_i\|_0.\quad\label{eq41} 
\end{eqnarray}

\noindent Again, definition \eqref{stiffmat}, property \eqref{vemstab-se2}, $\delta t<1$ and Cauchy-Schwarz inequality implies
\begin{eqnarray}
I_2:=(\delta t)\,a_h(\psi_i^n,L^\psi_i) \leq (\max\limits_{E}\Bar{\xi_i}^E)\,\max\lbrace1,\beta^*\rbrace\,\|\nabla\psi_i^n\|_0\,\|\nabla L^\psi_i\|_0.\quad\label{eq43} 
\end{eqnarray}
Let $\omega_{max}:=\max\lbrace\|\omega_1\|_\infty,\|\omega_2\|_\infty\rbrace$, using $\delta t<1$ and Cauchy-Schwarz inequality we get
\begin{eqnarray}
I_3&\leq&\omega_{max}\|\boldsymbol{\mathcal{P}}\nabla \psi^n_i\|_0\,\|\mathcal{P}L^\psi_i\|_0\leq \omega_{max}\|\nabla \psi^n_i\|_0\,\|L^\psi_i\|_0.\quad\label{eq45}
\end{eqnarray}

\noindent Adding and subtracting $(\delta t)\Big( \mathcal{P}W^n_i\big(\textstyle\sum\limits_{j=1}^m\mathsf{A}(i,j)\mathcal{P}Z^n_j\big),\,\mathcal{P}L^\psi_i\, \Big)_{0}$ to $I_4,$ we get
\begin{eqnarray}
I_4&=&(\delta t)\Big[\,\Big( \mathcal{P}W^n_i\big(\textstyle\sum\limits_{j=1}^m\mathsf{A}(i,j)\mathcal{P}(W^n_j-Z_j^n)\big),\,\mathcal{P}L^\psi_i\, \Big)_{0}+\Big( \mathcal{P}(W_i^n-Z^n_i)\big(\textstyle\sum\limits_{j=1}^m\mathsf{A}(i,j)\mathcal{P}Z^n_j\big),\,\mathcal{P}L^\psi_i\, \Big)_{0}\,\Big].\nonumber
\end{eqnarray}
Considering $\mathsf{A}_{max}:=\max\limits_{i,j}|\mathsf{A}(i,j)|$, using generalized H$\Ddot{o}$lder's inequality, $\delta t<1$ and triangle inequality we obtain
\begin{eqnarray}
I_4&\leq& \mathsf{A}_{max}\,\Big[\,\|\mathcal{P}W_i^n\|_{L^4(\Omega)}\,\textstyle\sum\limits_{j=1}^m\| \mathcal{P}\psi_j^n\|_{L^4(\Omega)} +\|\mathcal{P}\psi_i^n\|_{L^4(\Omega)}\,\textstyle\sum\limits_{j=1}^m\| \mathcal{P}Z_j^n\|_{L^4(\Omega)}\,\Big]\,\|\mathcal{P}L^\psi_i\|_{0}.\nonumber
\end{eqnarray}
Using \eqref{sobolevimbedding}, property of $\mathcal{P}$, Poincar$\Acute{e}$ inequality (with constant $C_P$), H$\Ddot{o}$lder's inequality and noting $\boldsymbol{\|Z\|}<\nu +\boldsymbol{\|W\|}$, we get
\begin{eqnarray}
I_4&\leq& C\,\mathsf{A}_{max}\,\Big[\,\|W_i^n\|_{1}\,\textstyle\sum\limits_{j=1}^m\|\psi_j^n\|_{1} +\|\psi_i^n\|_{1}\,\textstyle\sum\limits_{j=1}^m\|Z_j^n\|_{1}\,\Big]\,\|L^\psi_i\|_{0}\nonumber\\
&\leq& C\,\mathsf{A}_{max}\,\big[\,2\boldsymbol{\|W\|}+\nu\,\big]\boldsymbol{\|\psi\|\,\|L^\psi\|}.\quad\label{eq44}
\end{eqnarray}

Note that, by adding and subtracting $\mathcal{P}W_j^n\,\mathcal{P}Z_\ell^n$ we have
\begin{eqnarray}
[\,\mathcal{P}W_\ell^n\,\mathcal{P}W_j^n-\mathcal{P}Z_\ell^n\,\mathcal{P}Z_j^n\,]= \mathcal{P}\psi_\ell^n\,\mathcal{P}W_j^n+\mathcal{P}Z_\ell^n\,\mathcal{P}\psi_j^n.\qquad\label{e45}
\end{eqnarray}

Next, let $\mathsf{Q}_{max}:=\max\limits_{i,\ell,j}|\mathsf{Q}(i,\ell,j)|$, using $\delta t<1$, \eqref{e45} and similar to \eqref{eq44}, we have
\begin{eqnarray}
I_5&\leq&\textstyle  \mathsf{Q}_{max} \Big(\,\big[\sum\limits_{\ell\neq i} \mathcal{P}\psi^n_\ell\sum\limits_{j\neq i}\mathcal{P}W^n_j+\sum\limits_{\ell\neq i}\mathcal{P}Z^n_\ell\sum\limits_{j\neq i} \mathcal{P}\psi^n_j\,\big],\,\mathcal{P}L^\psi_i\, \Big)_{0}\nonumber\\
&\leq& C\,\mathsf{Q}_{max}\,\big[\,2\boldsymbol{\|W\|}+\nu\,\big]\boldsymbol{\|\psi\|\,\|L^\psi\|}.\quad\label{eq46}
\end{eqnarray}
Denoting $\mathsf{R}_{max}:=\max\limits_{i,j}|\mathsf{R}(i,j)|$, using $\delta t<1$, Cauchy-Schwarz inequality, triangle inequality and Poincar$\Acute{e}$ inequality, we obtain
\begin{eqnarray}
I_6&\leq&\textstyle \mathsf{R}_{max}\,\big(\sum\limits_{j=1}^m\|\mathcal{P}\psi_j^n\|_0\big)\,\|\mathcal{P}L^\psi_i\|_0\leq C_P^2\, \mathsf{R}_{max}\,\boldsymbol{\|\psi\|\,\|L^\psi\|}.\label{eq47}
\end{eqnarray}
Substituting the results  \eqref{eq41}-\eqref{eq45}, \eqref{eq44}, \eqref{eq46}-\eqref{eq47} into \eqref{eq4}, and using Poincar$\Acute{e}$ inequality, we obtain
\begin{eqnarray}
I&\leq&\big[\, C^2_P\,\max\lbrace1,\alpha^*\rbrace+(\max\limits_{E}\Bar{\xi_i}^E)\,\max\lbrace1,\beta^*\rbrace\,+C_P\,\omega_{max}\,\big]\,\|\nabla\psi_i^n\|_0\,\|\nabla L^\psi_i\|_0\nonumber\\
&\,&+\big[\,C\,\big(\,\mathsf{A}_{max}+\mathsf{Q}_{max}\,\big)\,\big(\,2\boldsymbol{\|W\|}+\nu\,\big)+C_P^2\,\mathsf{R}_{max}\,\big]\,\boldsymbol{\|\psi\|\,\|L^\psi\|}.\label{eq11}
\end{eqnarray}
Summing \eqref{eq11} over $i=1$ to $m$, and using H$\Ddot{o}$lder's inequality, we get \begin{eqnarray}
\boldsymbol{\|L(W)-L(Z)\|}&\leq& (\,\mathcal{C}_1\,\nu+\mathcal{C}_2\,)\boldsymbol{\|W-Z\|}, \label{eqcon}
\end{eqnarray}
 \begin{eqnarray}
 \text{where,}\qquad\mathcal{C}_1&=&m\,C\,\big(\,\mathsf{A}_{max}+\mathsf{Q}_{max}\,\big).\qquad\qquad\qquad
 \label{c1}\\
 \mathcal{C}_2&=&\big[\, C^2_P\,\max\lbrace1,\alpha^*\rbrace+(\max\limits_{E}\Bar{\xi_i}^E)\,\max\lbrace1,\beta^*\rbrace\,+C_P\,\omega_{max}\,\big]\nonumber\\
 &&\quad+m\,\big[\,C\,\big(\,\mathsf{A}_{max}+\mathsf{Q}_{max}\,\big)\,2\boldsymbol{\|W\|}+C_P^2\,\mathsf{R}_{max}\,\big].\label{c2}
 \end{eqnarray}
 Under suitable choice for $\nu$ (see Remark \ref{rem:continuity}), the estimate \eqref{eqcon} implies $\boldsymbol{\|L(W)-L(Z)\|}<\epsilon$. Hence $\boldsymbol{L}$ is continuous.
 
 \begin{eqnarray}
 \text{Consider the fixed constant,}\qquad\quad K=\dfrac{2\,\big(\,\max\lbrace 1,\alpha^*\rbrace\,\boldsymbol{\|U^{n-1}\|}+1\,\big)}{C\,\min\lbrace 1,\alpha_*\rbrace+(\delta t)\,\min\limits_{i}\lbrace 1,\beta_*, \xi_{i,0}\rbrace}.\qquad\qquad\qquad\label{kval}
 \end{eqnarray}
 
 We will show $\boldsymbol{\langle L(W),W\rangle}>0$, for all $\boldsymbol{W}\in \boldsymbol{\mathcal{V}_h^p}$ with $\boldsymbol{\|W\|}=K$. We have,$  \boldsymbol{\langle L(W),W\rangle}=\textstyle\sum\limits_{i=1}^m\mathcal{F}_i(\boldsymbol{W},\boldsymbol{W})$. Consider
 
 \begin{eqnarray}
\mathcal{F}_i(\boldsymbol{W},\boldsymbol{W})&=&m_h(W_i^n,W_i^n)-m_h(U_i^{n-1},W_i^n) +(\delta t)\, a_h(W_i^n,W_i^n)\nonumber\\
&&\quad+(\delta t)\Big( (\omega_1,\omega_2)\cdot\boldsymbol{\mathcal{P}}\nabla W^n_i,\,\mathcal{P}W_i^n\, \Big)_{0}+(\delta t)\Big( \mathcal{P}W^n_i\big(\textstyle\sum\limits_{j=1}^m\mathsf{A}(i,j)\mathcal{P}W^n_j\big),\,\mathcal{P}W_i^n\, \Big)_{0} \nonumber\\
&&\quad+(\delta t)\Big( \big(\textstyle\sum\limits_{\ell,j\neq i}\mathsf{Q}(i,\ell,j)\mathcal{P}W^n_\ell \mathcal{P}W^n_j\big),\,\mathcal{P}W_i^n\, \Big)_{0}+(\delta t)\Big( \big(\textstyle\sum\limits_{j=1}^m\mathsf{R}(i,j)\mathcal{P}W^n_j\big),\,\mathcal{P}W_i^n\, \Big)_{0} \nonumber\\
&&\quad- (\delta t)(\,f^n_i,\,\mathcal{P}W_i\,)_0 := \textstyle\sum\limits_{r=1}^8II_r.\quad\label{eq2e1}
 \end{eqnarray}
 
 Using Cauchy-Schwarz inequality, \eqref{vemstab-se1}, \eqref{vemstab-se2} and Poincar$\Acute{e}$ inequality, we  have
 \begin{eqnarray}
 \textstyle\sum\limits_{r=1}^3II_r\geq \min\lbrace 1,\alpha_*\rbrace\,\|W_i^n\|^2_0-\max\lbrace 1,\alpha^*\rbrace\,C_P^2\,\|\nabla W_i^n\|_0\,\|\nabla U^{n-1}_i\|_0+(\delta t)\,\xi_{i,0}\,\min\lbrace 1,\beta_*\rbrace\,\|\nabla W_i^n\|^2_0,\qquad\label{eq2e11}
 \end{eqnarray}
 where $\xi_{i,0}:=\min\limits_{x\in\Omega}|\xi_i(t_n,x)|.\,$
 Again using Cauchy-Schwarz inequality and Poincar$\Acute{e}$ inequality, we get
 \begin{eqnarray}
 |II_4|\leq (\delta t)\omega_{max}\,C_p\,\|\nabla W_i^n\|^2_0
 \end{eqnarray}
 Similar to \eqref{eq44}, \eqref{eq46} and \eqref{eq47}, we obtain,
 \begin{eqnarray}
 \textstyle\sum\limits_{r=5}^7|II_r|\leq (\delta t)\,C\,\big(\,\mathsf{A}_{max}+\mathsf{Q}_{max}\,\big)\,\boldsymbol{\|W\|}^3+(\delta t)\,\mathsf{R}_{max}\,C_P^2\,\boldsymbol{\|W\|}^2.\label{eq2e12}
 \end{eqnarray}
 Using Cauchy-Schwarz inequality and Poincar$\Acute{e}$ inequality, we get
 \begin{eqnarray}
 II_8\geq -(\delta t)\,C_P\,\|f_i^n\|_0\,\|\nabla W_i^n\|_0.\label{eq2e13}
 \end{eqnarray}
 
 From the estimates \eqref{eq2e11}-\eqref{eq2e13}, we obtain,
 \begin{eqnarray}
 &&\mathcal{F}_i(\boldsymbol{W},\boldsymbol{W})\geq \min\lbrace 1,\alpha_*\rbrace\,\|W_i^n\|^2_0-\max\lbrace 1,\alpha^*\rbrace\,\|W_i^n\|_0\,\|U^{n-1}_i\|_0+(\delta t)\,\xi_{i,0}\,\min\lbrace 1,\beta_*\rbrace\,\|\nabla W_i^n\|^2_0\nonumber\\
 &&\quad- (\delta t)\omega_{max}\,C_p\,\|\nabla W_i^n\|^2_0-
 (\delta t)\,C\,\big(\,\mathsf{A}_{max}+\mathsf{Q}_{max}\,\big)\,\boldsymbol{\|W\|}^3-(\delta t)\,\mathsf{R}_{max}\,C_P^2\,\boldsymbol{\|W\|}^2 \nonumber\\
 &&\quad-(\delta t)\,C_P\,\|f_i^n\|_0\,\|\nabla W_i^n\|_0.\nonumber
 \end{eqnarray}
 Summing $\mathcal{F}_i(\boldsymbol{W},\boldsymbol{W})$ over $i=1$ upto $m$, equivalance of $\|\cdot\|_0$, $\|\nabla\cdot\|_0$ in $V_h^p$ and using H$\Ddot{o}$lder's inequality, we get
 \begin{eqnarray}
 &&\boldsymbol{\langle L(W),W\rangle}\geq \boldsymbol{\|W\|}\Big\lbrace\,\big(\, C\,\min\lbrace 1,\alpha_*\rbrace+(\delta t)\,\min\limits_{i}\lbrace 1,\beta_*, \xi_{i,0}\rbrace\,\big)\,\boldsymbol{\|W\|}-\max\lbrace 1,\alpha^*\rbrace\,\boldsymbol{\|U^{n-1}\|}\nonumber\\
 &&\quad-(\delta t)\big[\,\omega_{max}\,C_p\,\boldsymbol{\|W\|}+C\,\big(\,\mathsf{A}_{max}+\mathsf{Q}_{max}\,\big)\,\boldsymbol{\|W\|}^2+\mathsf{R}_{max}\,C_P^2\,\boldsymbol{\|W\|}+C_P\,\boldsymbol{\|f\|}\,\big]\,\Big\rbrace.\label{eqsign}
 \end{eqnarray}
 For sufficiently small $(\delta t)$,  we can have the estimate
 \begin{eqnarray}
 (\delta t)\big[\,\omega_{max}\,C_p\,K+C\,\big(\,\mathsf{A}_{max}+\mathsf{Q}_{max}\,\big)\,K^2+\mathsf{R}_{max}\,C_P^2\,K+C_P\,\boldsymbol{\|f\|}\,\big] <1.\label{amp1}
 \end{eqnarray}
 Substituting \eqref{amp1} into \eqref{eqsign} and using \eqref{kval}, we obtain
 
 \begin{eqnarray}
 %\boldsymbol{\langle L(W),W\rangle}&\geq& K\,\Big\lbrace\,\Big(\, C\,\min\lbrace 1,\alpha_*\rbrace+(\delta t)\,\min\limits_{i}\lbrace 1,\beta_*, \xi_{i,0}\rbrace\,\Big)\,K-\max\lbrace 1,\alpha^*\rbrace\,\boldsymbol{\|U^{n-1}\|}-1\,\Big\rbrace\nonumber\\
\boldsymbol{\langle L(W),W\rangle} &>&0\quad\text{for any $\boldsymbol{W}\in\boldsymbol{\mathcal{V}_h^p}$ with }\,\boldsymbol{\|W\|}=K.\nonumber
 \end{eqnarray}
 Thus by Proposition \ref{prop-Brouwer}, existence of a solution $\boldsymbol{U^n}\in\boldsymbol{\mathcal{V}_h^p}$ to the discrete problem \eqref{vemform} is proved.
 
 Next we show uniqueness of the solution $\boldsymbol{U^n}\in\boldsymbol{\mathcal{V}_h^p}$ to \eqref{vemform} by the method of contradiction. Let $\boldsymbol{U^n},\,\boldsymbol{\widetilde{U}^n}\in\boldsymbol{\mathcal{V}_h^p}$ be two distinct solutions of \eqref{vemform}. Denote $\boldsymbol{E^n}:=\boldsymbol{U^n}-\boldsymbol{\widetilde{U}^n}\in\boldsymbol{\mathcal{V}_h^p}$. Then we have 
 \begin{eqnarray}
 0 &=&\sum\limits_{i=1}\bigg\lbrace m_h(E_i^n,E^n_i)+(\delta t)\, a_h(E^n_i,E^n_i)+(\delta t)\Big( (\omega_1,\omega_2)\cdot\boldsymbol{\mathcal{P}}\nabla E^n_i,\,\mathcal{P}E^n_i\, \Big)_{0}\nonumber\\
&&\quad+(\delta t)\Big( \mathcal{P}U^n_i\big(\textstyle\sum\limits_{j=1}^m\mathsf{A}(i,j)\mathcal{P}U^n_j\big)-\mathcal{P}\widetilde{U}^n_i\big(\textstyle\sum\limits_{j=1}^m\mathsf{A}(i,j)\mathcal{P}\widetilde{U}^n_j\big),\,\mathcal{P}E^n_i\, \Big)_{0} \nonumber\\
&&\quad+(\delta t)\Big( \big(\textstyle\sum\limits_{\ell,j\neq i}\mathsf{Q}(i,\ell,j)\big[\,\mathcal{P}U^n_\ell \mathcal{P}U^n_j-\mathcal{P}\widetilde{U}^n_\ell \mathcal{P}\widetilde{U}^n_j\,\big]\,\big),\,\mathcal{P}E^n_i\, \Big)_{0}+(\delta t)\Big( \big(\textstyle\sum\limits_{j=1}^m\mathsf{R}(i,j)\,\mathcal{P}E^n_j\,\big),\,\mathcal{P}E^n_i\, \Big)_{0}\,\bigg\rbrace.\nonumber
 \end{eqnarray}
 From previous estimates we have
 \begin{eqnarray}
 0&\geq& \big[\, C\,\min\lbrace 1,\alpha_*\rbrace+(\delta t)\,\min\limits_{i}\lbrace 1,\beta_*, \xi_{i,0}\rbrace\,-(\delta t)\big[\,\omega_{max}\,C_p\nonumber\\
 &&\quad+(\delta t)\,C\,\big(\,\mathsf{A}_{max}+\mathsf{Q}_{max}\,\big)\,(\boldsymbol{\|U^n\|}+\boldsymbol{\|\widetilde{U}^n\|})+(\delta t)\,\mathsf{R}_{max}\,C_P^2\,\big]\,\boldsymbol{\|E^n\|}^2.\nonumber
 \end{eqnarray}
 Choose sufficiently small $(\delta t)$ such that
 \begin{eqnarray}
C\,\min\lbrace 1,\alpha_*\rbrace+(\delta t)\,\min\limits_{i}\lbrace 1,\beta_*, \xi_{i,0}\rbrace&\geq&(\delta t)\big[\,\omega_{max}\,C_p+C\,\big(\,\mathsf{A}_{max}+\mathsf{Q}_{max}\,\big)\,(\boldsymbol{\|U^n\|}+\boldsymbol{\|\widetilde{U}^n\|})\nonumber\\
 &&\qquad\qquad\qquad+\mathsf{R}_{max}\,C_P^2\,\big]\nonumber
 \end{eqnarray}
  implies $\boldsymbol{\|E^n\|}^2\leq 0$. Therefore $\boldsymbol{U^n}=\boldsymbol{\widetilde{U}^n}$.
\end{proof}

\begin{Rem}
Let any $\epsilon>0$ be given and note that $\mathcal{C}_1>0,\,\,\mathcal{C}_2>0$. Choose $0<\epsilon_*<\epsilon$. Then note that the quadratic equation $\mathcal{C}_1\,\nu^2+\mathcal{C}_2\,\nu-\epsilon_*=0$ has real positive solution. This guarantees a $\nu>0$ satisfying $\nu\leq \dfrac{\epsilon}{\mathcal{C}_1\,\nu+\mathcal{C}_2}$.\label{rem:continuity} 
\end{Rem}

\section{Linearized scheme}
\label{sec:linearscheme}
Considering the implementation aspects, solving highly nonlinear systems, especially when the systems are time-dependent, is cumbersome and time-consuming. Hence we modify the nonlinear scheme \eqref{vemform} and an equivalent  linear formulation which is stated as follows :

\noindent for each $n=1,2,...,N,$\,\,\, Find\,\, $(U_1^n,U_2^n,...,U_m^n)\in \boldsymbol{\mathcal{V}_h^p}$\,\,\, such that \quad  $\forall\,\boldsymbol{V}\in \boldsymbol{\mathcal{V}_h^p}\, :$
\begin{eqnarray}
&& U^0_i(x)\, =\, u_i^0(x) \qquad\quad \text{in}\,\,\,\Omega,\,\quad i=1,...,m,\nonumber\\
&&\textstyle\sum\limits_{i=1}^m\Big\lbrace\,m_h(U_i^n,V_i)-m_h(U_i^{n-1},V_i)+(\delta t)\,a_h(U_i^n,V_i)+(\delta t)\,l_{1,h}(U_i^n,V_i)+(\delta t)\,l_{2,h}(U_i^n,V_i)\nonumber\\
&&\quad+(\delta t)
\,l_{3,h}(U_i^n,V_i)\Big\rbrace \, =\, \textstyle (\delta t)\sum\limits_{i=1}^m\Big\lbrace-
l_{4,h}(U_i^n,V_i)- l_{5,h}(U_i^n,V_i)+ (\,f^n_i,\,\mathcal{P}V_i\,)_0\,\Big\rbrace, \label{vemformlinear}
\end{eqnarray}
where,

\begin{eqnarray}
&&l_{1,h}(U_i^n,V_i) =\Big( (\omega_1,\omega_2)\cdot\boldsymbol{\mathcal{P}}\nabla U^n_i,\,\mathcal{P}V_i\, \Big)_{0},\quad
l_{2,h}(U_i^n,V_i)=\Big( \mathcal{P}U^n_i\big(\textstyle\sum\limits_{j=1}^m\mathsf{A}(i,j)\mathcal{P}U^{n-1}_j\big),\,\mathcal{P}V_i\, \Big)_{0},\nonumber\\
&&l_{3,h}(U_i^n,V_i)=\Big(\mathsf{R}(i,i)\mathcal{P}U^n_i,\,\mathcal{P}V_i\, \Big)_{0},\quad l_{4,h}(U_i^n,V_i)=\Big( \big(\textstyle\sum\limits_{\ell,j\neq i}\mathsf{Q}(i,\ell,j)\mathcal{P}U^{n-1}_\ell \mathcal{P}U^{n-1}_j\big),\,\mathcal{P}V_i\, \Big)_{0},\nonumber\\
&&l_{5,h}(U_i^n,V_i)=\Big( \big(\textstyle\sum\limits_{j\neq i}\mathsf{R}(i,j)\mathcal{P}U^{n-1}_j\big),\,\mathcal{P}V_i\, \Big)_{0}.\, \nonumber
\end{eqnarray}

Consider a finite sequence $\lbrace \boldsymbol{\phi^n}\rbrace_{n=1}^N$ of functions in $\boldsymbol{\mathcal{V}_h^p}$ associated to different time levels. In our analysis we use the norm $\boldsymbol{[\phi]}_{0,k}$ that is defined as,
\begin{eqnarray}
\textstyle
\boldsymbol{[\phi]}_{0,k}:=\Big(\,(\delta t)  \sum\limits_{n=1}^N\Big(\,\sum\limits_{i=1}^m\|\phi_i^n\|^2_k\,\Big)\,\Big)^\frac{1}{2}.
\end{eqnarray}

In this section, we shall first discuss the well-posedness of the linear virtual element formulation \eqref{vemformlinear}. Then we derive a priori  error estimates involving the numerical solution of  \eqref{vemformlinear} with respect to the norm  $\boldsymbol{[\cdot]}_{0,1}$.

\begin{Th}
Under the assumption \eqref{assump1} and 
for sufficiently small $(\delta t)$, we assume there exists a unique $\boldsymbol{U^k}\in\boldsymbol{\mathcal{V}_h^p}$ for $0\leq k\leq n-1$ satisfying  \eqref{vemformlinear}. Then for all $1\leq n\leq N$ there exists a unique solution $\boldsymbol{U^n} \in\boldsymbol{\mathcal{V}_h^p}$ to the linear discrete problem \eqref{vemformlinear}.
\end{Th}
\begin{proof}
The proof is analogous to the proof of Theorem \ref{thm:nonlinear}. 
\end{proof}

The following are some valid results that will help prove our error estimate.  
Consider the local polynomial interpolation estimates (see Lemma 5.1 in\,\cite{vem4}) : for all $ E\in\mathcal{T}_h$ and any $w\in H^{s+1}(E)$,

\begin{eqnarray}
	&\,&\|w-\mathcal{P}w\|_{m,E}\leq C\,h_E^{s+1-m}\,|w|_{s+1,E}\quad m,s\in\mathbb{N}\cup \lbrace 0\rbrace,\,\,m\leq s+1\leq k+1.\label{pi0interpol}\\
	&\,&\|w-\mathcal{P}_k^\nabla w\|_{m,E}\leq C\,h_E^{s+1-m}\,|w|_{s+1,E}\quad m,s\in\mathbb{N},\,\,m\leq s+1\leq k+1,\,\,s\geq 1. \label{pigradinterpol}
\end{eqnarray}

\noindent The virtual interpolation estimate given below can be found in \cite{veminfinity}. For $1\leq s\leq k,\,\,\forall\,E\in\mathcal{T}_h\,$ and for every $w\in H^{1+s}(E)$, there exists $\mathcal{W}\in V_h^k$ satisfying
\begin{eqnarray}
	\|w-\mathcal{W}\|_{E}+h\,|w-\mathcal{W}|_{1,E}+h\,\|w-\mathcal{W}\|_{\infty,E}\leq C\,h^{1+s}\,|w|_{1+s,\Omega}.\label{veminterpol}
\end{eqnarray} 

\begin{Rem}
On a bounded domain $E$ with area of $E\,\approx\,h_{E}^2,$ and also if $\|v\|_{\infty}<\infty,\,\|v\|_{0}<\infty$, then using inverse inequality on polynomials, property of $\mathcal{P}$ and standard estimation, we have
\begin{eqnarray}
\|\mathcal{P}v\|_{\infty}\leq h_{E}^{-1}\|\mathcal{P}v\|_{0}\leq h_{E}^{-1}\|v\|_{0}\leq h_{E}^{-1}h_{E}\|v\|_{\infty}\leq\|v\|_{\infty}.\quad\label{projlinf} 
\end{eqnarray}
\end{Rem}

\begin{Th}
Let $\boldsymbol{u^n}\in [ H^{s+1}(\Omega)]^m,\,\,s\geq 1$ be the smooth exact solution satisfying \eqref{modelprob}-\eqref{ic} and $\boldsymbol{U^n}\in\boldsymbol{\mathcal{V}_h^p}$ be the solution of discrete form \eqref{vemformlinear} at $n^{th}$ time step. Under the assumptions \eqref{assump} and \eqref{assump1}, the numerical solution $\boldsymbol{U^n}$ converges to the exact solution $\boldsymbol{u^n}$, as $h\rightarrow 0$ and we obtain the estimates,
\begin{eqnarray}
%\boldsymbol{[U^n-u^n]}_{\infty,0}&\leq&A_1h^s+A_2(\delta t).\label{erres1}%\\
\boldsymbol{[U^n-u^n]}_{0,1}&\leq&A_1h^s+A_2(\delta t).\label{erres2}
\end{eqnarray}
where $A_1$ and $A_2$ are independent of $h$ and $\delta t.$
\end{Th}
\begin{proof}
Let $\boldsymbol{\mathcal{U}^n}\in \boldsymbol{\mathcal{V}_h^p}$ be the interpolant of $\boldsymbol{u^n}$ onto $\boldsymbol{\mathcal{V}_h^p}$. We denote,
\begin{eqnarray}
\boldsymbol{\varphi^n}:=\boldsymbol{u^n-U^n}=(\,\boldsymbol{u^n-\mathcal{U}^n}\,)+(\,\boldsymbol{\mathcal{U}^n-U^n}\,)=\boldsymbol{\eta^n}+\boldsymbol{\chi^n}.\label{splite}
\end{eqnarray}
We need to estimate $\boldsymbol{\chi^n}\in \boldsymbol{\mathcal{V}_h^p}$. Since $\boldsymbol{u^n}$ also satisfies the variational form \eqref{weakform}, $\boldsymbol{U^n}$ satisfies \eqref{vemformlinear} and  rewriting, we have\,\,\,$ \forall\,\boldsymbol{V}\in \boldsymbol{\mathcal{V}_h^p}$,
\begin{eqnarray}
&&\textstyle\sum\limits_{i=1}^m\bigg\lbrace\,\underbrace{\left(\frac{\partial u^n_i}{\partial t}, V_i\right)_{0}-\frac{1}{\delta t}\,m_h(U^n_i-U^{n-1}_i,V_i)}_{:=I}\,+\,\underbrace{\left(\, \xi_i\nabla u^n_i,\nabla V_i\,\right)_{0}-a_h(U_i^n,V_i)}_{:=II}\,\nonumber\\
&&\quad+\left( (\omega_1,\omega_2)\cdot\nabla u^n_i,\,V_i\, \right)_{0}-l_{1,h}(U_i^n,V_i)\,+\,\Big( u^n_i\big(\textstyle\sum\limits_{j=1}^m\mathsf{A}(i,j)u^n_j\big),\,V_i\, \Big)_{0}-l_{2,h}(U_i^n,V_i)\,\nonumber\\
&&\quad + \,\left(\mathsf{R}(i,i)u_i^n,\,V_i\, \right)_{0}    - l_{3,h}(U_i^n,V_i)\,=\, -\,\Big( \big(\textstyle\sum\limits_{\ell,j\neq i}\mathsf{Q}(i,\ell,j)u^n_\ell u^n_j\big),\,V_i\, \Big)_{0}-l_{4,h}(U_i^n,V_i)\,\nonumber\\
&&\quad\quad\quad\, -\,\Big( \big(\textstyle\sum\limits_{j\neq i}^m\mathsf{R}(i,j)u^n_j\big),\,V_i\, \Big)_{0}-l_{5,h}(U_i^n,V_i)\,+ \left(\,f^n_i,\,(I-\mathcal{P})V_i\,\right)_0\,\bigg\rbrace.\qquad\label{er-e1es}
\end{eqnarray}

 Adding and subtracting the term $\dfrac{1}{\delta t}\,m_h(u^n_i-u^{n-1}_i,V_i)$ to $I$ and using \eqref{splite}, we get,
 \begin{eqnarray}
 I=\left(\frac{\partial u^n_i}{\partial t}, V_i\right)_{0}-\frac{1}{\delta t}\,m_h(u^n_i-u^{n-1}_i,V_i)+\frac{1}{\delta t}\,m_h(\eta^n_i-\eta^{n-1}_i,V_i)+\frac{1}{\delta t}\,m_h(\chi^n_i-\chi^{n-1}_i,V_i).\qquad\label{er-e1es1}
 \end{eqnarray}
 Next to $II$, adding and subtracting $a_h(u_i^n,V_i)$ and using \eqref{splite}, we obtain
 \begin{eqnarray}
 II=\left(\, \xi_i\nabla u^n_i,\nabla V_i\,\right)_{0}-a_h(u_i^n,V_i)+a_h(\eta_i^n,V_i)+a_h(\chi_i^n,V_i).\qquad\label{er-e1es2}
 \end{eqnarray}
 Substituting \eqref{er-e1es1}-\eqref{er-e1es2} into \eqref{er-e1es}, letting $V_i:=\chi_i^n$ and rearranging, we have
 
 \begin{eqnarray}
&&\textstyle\sum\limits_{i=1}^m\bigg\lbrace\,m_h(\chi^n_i-\chi^{n-1}_i,\chi_i^n)+\delta t\,a_h(\chi_i^n,\chi^n_i) \nonumber\\
&&\quad=\,-\delta t\,\left(\frac{\partial u^n_i}{\partial t}, \chi_i^n\right)_{0}+m_h(u^n_i-u^{n-1}_i,\chi_i^n)-m_h(\eta^n_i-\eta^{n-1}_i,\chi_i^n)-\delta t\,a_h(\eta_i^n,\chi^n_i)\nonumber\\
&&\quad\quad+\delta t\,a_h(u_i^n,\chi^n_i)-\delta t\,\left(\, \xi_i\nabla u^n_i,\nabla \chi^n_i\,\right)_{0}+\delta t\,\Big[\,l_{1,h}(U_i^n,\chi_i^n)-\left( (\omega_1,\omega_2)\cdot\nabla u^n_i,\,\chi_i^n\, \right)_{0}\,\Big]\nonumber\\
&&\quad\quad-\delta t\,\Big[\,\Big( u^n_i\big(\textstyle\sum\limits_{j=1}^m\mathsf{A}(i,j)u^n_j\big),\,\chi_i^n\, \Big)_{0}-l_{2,h}(U_i^n,\chi_i^n)\,\Big] -\delta t\, \Big[\,\left(\mathsf{R}(i,i)u_i^n,\,\chi_i^n\, \right)_{0}    - l_{3,h}(U_i^n,\chi_i^n)\,\Big]   \nonumber\\
&&\quad\quad-\delta t\,\Big[\,\Big( \big(\textstyle\sum\limits_{\ell,j\neq i}\mathsf{Q}(i,\ell,j)u^n_\ell u^n_j\big),\,\chi_i^n\, \Big)_{0}-l_{4,h}(U_i^n,\chi_i^n)\,\Big]\nonumber\\
&&\quad\quad\, -\delta t\,\Big[\,\Big( \big(\textstyle\sum\limits_{j\neq i}^m\mathsf{R}(i,j)u^n_j\big),\,\chi_i^n\, \Big)_{0}-l_{5,h}(U_i^n,\chi_i^n)\,\Big]+\,\delta t\, \left(\,f^n_i,\,(I-\mathcal{P})\chi_i^n\,\right)_0\,\bigg\rbrace. \qquad \label{mx1}
\end{eqnarray}

%Using \eqref{vemstab-se1} and the identity $\frac{1}{2}(a^2-b^2)\leq a(a-b)$, we obtain,
%\begin{eqnarray}
%m_h(\chi^n_i-\chi^{n-1}_i,\chi_i^n)\geq \min\lbrace1,\alpha_*\rbrace\,(\,\chi^n_i-\chi^{n-1}_i,\chi_i^n\,)_0\geq \frac{\min\lbrace1,\alpha_*\rbrace}{2}\,\big[\,\|\chi_i^n\|_0^2- \|\chi_i^{n-1}\|_0^2\,\big].\qquad\qquad\label{x3}
%\end{eqnarray}
%%%%%%%%%%%%%%%%%
\noindent Using \eqref{vemstab-se1}, property of $\mathcal{P}$ and   $mn\leq \dfrac{m^2}{\theta}+n^2\theta\,\,(\text{with}\,\theta:={\min\lbrace1,\alpha_*\rbrace}/{2}\,)$ we have

\begin{eqnarray}
m_h(\chi_i^n,\chi_i^n)&\geq& \min\lbrace1,\alpha_*\rbrace\,\|\chi_i^n\|_0^2.\quad\label{x1}\\
m_h(\chi_i^{n-1},\chi_i^n)&\leq& (1+\alpha^*)\,\|\chi_i^{n-1}\|_0\,\|\chi_i^n\|_0\leq \frac{2(1+\alpha^*)^2}{\min\lbrace1,\alpha_*\rbrace}\,\|\chi_i^{n-1}\|_0^2+\frac{\min\lbrace1,\alpha_*\rbrace}{2}\|\chi_i^n\|_0^2.\qquad\qquad\label{x2}
\end{eqnarray}
Then \eqref{x1} and \eqref{x2} implies 
\begin{eqnarray}
m_h(\chi^n_i-\chi^{n-1}_i,\chi_i^n)\geq \frac{\min\lbrace1,\alpha_*\rbrace}{2}\|\chi_i^n\|_0^2- \frac{2(1+\alpha^*)^2}{\min\lbrace1,\alpha_*\rbrace}\,\|\chi_i^{n-1}\|_0^2.\label{x3}
\end{eqnarray}
%%%%%%%%%%%%%%%%
Using definition \eqref{stiffmat} and \eqref{vemstab-se2}, we obtain
\begin{eqnarray}
\delta t\,a_h(\chi_i^n,\chi_i^n)\geq \delta t\,\min\lbrace 1,\beta_*\rbrace\,\xi_{i,0}\,\|\nabla \chi_i^n\|^2_0.\label{x4}
\end{eqnarray}
 Using \eqref{massmat}, add and subtract $(u_i^n-u_i^{n-1},\,\mathcal{P}\chi_i^n+\chi_i^n)_0,\,$ and property of $\mathcal{P}$,\, we get
 \begin{eqnarray}
  I_1&:=&m_h(u^n_i-u^{n-1}_i,\chi_i^n)-\delta t\,\left(\frac{\partial u^n_i}{\partial t}, \chi_i^n\right)_{0}\nonumber\\
  &=&S_1((I-\mathcal{P})(u_i^n-u_i^{n-1}),(I-\mathcal{P})\chi_i^n)+(u_i^n-u_i^{n-1},\,(I-\mathcal{P})\chi_i^n)_0\nonumber\\
  &&\quad+\Big[\,(u_i^n-u_i^{n-1},\,\chi_i^n)-\delta t\,\left(\frac{\partial u^n_i}{\partial t}, \chi_i^n\right)_{0}\,\Big].\nonumber
 \end{eqnarray}
 Using \eqref{vemstab-se1}, Cauchy-Schwarz inequality and Taylor's theorem w.r.t $t$ ( neglecting $(\delta t)^2$ term ), we get
 \begin{eqnarray}
 I_1&\leq&(\alpha^*+1)\,\|u_i^n-u_i^{n-1}\|_0\,\|\chi_i^n\|_0+\Big[\,(u_i^n-u_i^{n-1},\,\chi_i^n)-\delta t\,\left(\frac{\partial u^n_i}{\partial t}, \chi_i^n\right)_{0}\,\Big]\nonumber\\
 &\leq& (\alpha^*+1)\,\delta t\,\| (u_i)^{n}_{t} \|_0\, \|\chi_i^n\|_0+\delta t\,\Big(\,\frac{u_i^n-u_i^{n-1}}{\delta t}-\frac{\partial u^n_i}{\partial t},\, \chi_i^n\,\Big)_{0}.\label{eqt-i1}
 \end{eqnarray}
 By applying Taylor series for variable $t$ with the remainder in terms of integral and using Cauchy-Schwarz inequality for the second term of \eqref{eqt-i1}, we have
 \begin{eqnarray}
 I_1 &\leq& (\alpha^*+1)\,\delta t\,\| (u_i)^{n}_{t} \|_0\, \|\chi_i^n\|_0+(\delta t)^{\frac{3}{2}}\,\|(u_i)_{tt}\|_0\,\|\chi_i^n\|_0.\qquad \label{eq-it1}
 \end{eqnarray}
 Using the  inequality $mn\leq \dfrac{m^2}{\theta}+n^2\theta\,$\,\,\text{with}\, $\theta:={\min\lbrace1,\alpha_*\rbrace}/{8(\alpha^*+1)}\,$ for the first term and $\theta:=1$ for the second term of \eqref{eq-it1}, we get,
 \begin{eqnarray}
 I_1&\leq& (\delta t)^2\,\frac{8(\alpha^*+1)^2}{\min\lbrace\alpha^*,1\rbrace}\,\| (u_i)_{t}^{n} \|_0^2+\frac{\min\lbrace\alpha_*,1\rbrace}{8}\, \|\chi_i^n\|_0^2+(\delta t)^2\,\|(u_i)_{tt}\|_0^2+\delta t\,\|\chi_i^n\|_0^2.\nonumber\\
 &&\label{eq-i1}
 \end{eqnarray}
 
 Using \eqref{vemstab-se1}, \eqref{vemstab-se2}, Cauchy-Schwarz inequality and the  inequality $mn\leq \dfrac{m^2}{\theta}+n^2\theta\,$\,\,\text{with}\, $\theta:={\min\lbrace1,\alpha_*\rbrace}/{8(\alpha^*+1)}\,)$ for the first term and $\theta:=(\xi_{i,0}\,{\min\lbrace1,\beta_*\rbrace})/{8(\beta^*+1)}\,)$ for the second term, we get,
 
 \begin{eqnarray}
I_2&:=&|m_h(\eta^n_i-\eta^{n-1}_i,\chi_i^n)|+|\delta t\,a_h(\eta_i^n,\chi^n_i)|\nonumber\\
&\leq& (1+\alpha^*)\,\|\eta^n_i-\eta^{n-1}_i\|_0\, \|\chi_i^n\|_0+(\delta t)\,(1+\beta^*)\, \|\nabla\eta_i^n\|_0\,\|\nabla\chi^n_i\|_0\nonumber\\
&\leq& \frac{8(\alpha^*+1)^2}{\min\lbrace\alpha_*,1\rbrace}\,2\,\big(\,\| \eta^n_i\|_0^2+\|\eta^{n-1}_i \|_0^2\,\big)+ \frac{\min\lbrace\alpha_*,1\rbrace}{8}\, \|\chi_i^n\|_0^2+\delta t \,\frac{8\,(1+\beta^*)^2}{\xi_{i,0}\,\min\lbrace\,1,\beta_*\,\rbrace}\,\|\nabla \eta_i^n\|_0^2\nonumber\\
&&\,+\delta t \,\frac{\xi_{i,0}\,\min\lbrace\,1,\beta_*\,\rbrace}{8}\,\|\nabla\chi_i^n\|_0^2\nonumber\\
&\leq& \frac{8(\alpha^*+1)^2}{\min\lbrace\alpha_*,1\rbrace}\,2\,C\,h^{2s+2}\,\big(\,\| u^n_i\|_{s+1}^2+\|u^{n-1}_i \|_{s+1}^2\,\big)+ \frac{\min\lbrace\alpha_*,1\rbrace}{8}\, \|\chi_i^n\|_0^2\qquad \text{( use \eqref{veminterpol} )}\,\nonumber\\
&&\,+\delta t \,\frac{8\,(1+\beta^*)^2}{\xi_{i,0}\,\min\lbrace\,1,\beta_*\,\rbrace}\,C\,h^{2s}\|u_i^n\|_{s+1}^2+\delta t \,\frac{\xi_{i,0}\,\min\lbrace\,1,\beta_*\,\rbrace}{8}\,\|\nabla\chi_i^n\|_0^2. \qquad\label{eq-i2}
 \end{eqnarray}
 
 Adding and subtracting $\delta t\,\left(\, \xi_i\nabla u^n_i,\boldsymbol{\mathcal{P}} \nabla \chi^n_i\,\right)_{0}=\delta t\, \left(\,\boldsymbol{\mathcal{P}} \xi_i\nabla u^n_i,\nabla \chi^n_i\,\right)_{0}$ we get
 \begin{eqnarray}
 I_3&:=&\delta t\,a_h(u_i^n,\chi^n_i)-\delta t\,\left(\, \xi_i\nabla u^n_i,\nabla \chi^n_i\,\right)_{0}\nonumber\\
 &=& \delta t\,\left(\, \xi_i\boldsymbol{\mathcal{P}}\nabla u^n_i-\xi_i\nabla u^n_i,\boldsymbol{\mathcal{P}}\nabla \chi^n_i\,\right)_{0}+\delta t\,\left(\, \boldsymbol{\mathcal{P}}(\xi_i\nabla u^n_i)-\xi_i\nabla u^n_i,\nabla \chi^n_i\,\right)_{0}\nonumber\\
 &&\quad+\delta t\,S_2((I-\Pi_p^\nabla)u_i^n,(I-\Pi_p^\nabla)\chi_i^n)\nonumber
 \end{eqnarray}
 Using Cauchy-Schwarz inequality, \eqref{vemstab-se2},  \eqref{pi0interpol} , \eqref{pigradinterpol} and  the  inequality $mn\leq \dfrac{m^2}{\theta}+n^2\theta\,$\,\,\text{with}\, $\theta:=\frac{\xi_{i,0}\,\min\lbrace\,1,\beta_*\,\rbrace}{8}$, we get
 \begin{eqnarray}
 I_3&\leq&\delta t\,\max\limits_E\Bar{\xi}_i^E\,\|(\boldsymbol{\mathcal{P}}-I)\nabla u_i^n\|_0\,\|\boldsymbol{\mathcal{P}}\nabla\chi_i^n\|_0+\delta t\,\|\boldsymbol{\mathcal{P}}(\xi_i\nabla u^n_i)-\xi_i\nabla u^n_i\|_0\,\|\nabla \chi^n_i\|_{0}\nonumber\\
 &&\quad+\delta t\,\beta^*\,\|\nabla(I-\Pi_p^\nabla)u_i^n\|_0\,\|\nabla(I-\Pi_p^\nabla)\chi_i^n\|_0\nonumber\\
 &\leq& \delta t\,\Big\lbrace\max\limits_E\Bar{\xi}_i^E\,\|\nabla (\Pi_p^\nabla-I)u_i^n\|_0+C\,h^s\,|\xi_i\nabla u_i^n|_{s}+
 \,\beta^*\,\|\nabla(I-\Pi_p^\nabla)u_i^n\|_0\,\Big\rbrace\,\|\nabla \chi^n_i\|_{0}\nonumber\\
 &\leq& \delta t\,\big(\max\limits_E\Bar{\xi}_i^E+\|\frac{\partial^s\xi_i}{\partial x^s}\|_\infty+\beta^*\big)\,C\,h^s\,\|u_i^n\|_{s+1}\,\|\nabla \chi^n_i\|_{0}\nonumber\\
 &\leq& \big(\max\limits_E\Bar{\xi}_i^E+\|\frac{\partial^s\xi_i}{\partial x^s}\|_\infty+\beta^*\big)^2\,\delta t\,\frac{8}{\xi_{i,0}\,\min\lbrace\,1,\beta_*\,\rbrace}\,C\,h^{2s}\,\|u_i^n\|^2_{s+1}+ \delta t\,\frac{\xi_{i,0}\,\min\lbrace\,1,\beta_*\,\rbrace}{8}\,\|\nabla \chi^n_i\|^2_{0}.\nonumber\\
 &&\label{eq-i3}
 \end{eqnarray}
 
 Next adding and  subtracting the term $\delta t\,\left( (\omega_1,\omega_2)\cdot\boldsymbol{\mathcal{P}}\nabla u^n_i,\,\mathcal{P}\chi_i^n\, \right)_{0}$,using Lemma 1 of \cite{vem8}, Cauchy-Schwarz inequality and triangle inequality, we obtain
 \begin{eqnarray}
 I_4&:=&\delta t\,\left( (\omega_1,\omega_2)\cdot\boldsymbol{\mathcal{P}}\nabla U^n_i,\,\mathcal{P}\chi_i^n\, \right)_{0}-\delta t\,\left( (\omega_1,\omega_2)\cdot\nabla u^n_i,\,\chi_i^n\, \right)_{0}\,\nonumber\\
 &=&\delta t\,\left( (\omega_1,\omega_2)\cdot\boldsymbol{\mathcal{P}}\nabla u^n_i,\,\mathcal{P}\chi_i^n\, \right)_{0}-\delta t\,\left( (\omega_1,\omega_2)\cdot\nabla u^n_i,\,\chi_i^n\, \right)_{0}+\delta t\,\left( (\omega_1,\omega_2)\cdot\boldsymbol{\mathcal{P}}\nabla \varphi^n_i,\,\mathcal{P}\chi_i^n\, \right)_{0}\nonumber\\
 &\leq& C\,\delta t\,\|(\omega_1,\omega_2)\|_{\infty,s}\,h^{s+1}\,\|u_i^n\|_{s+1}\,\|\nabla \chi_i^n\|_0+\omega_{max}\,\delta t\,\|\nabla\eta_i^n\|_0\,\|\chi_i^n\|_0+\omega_{max}\,\delta t\,\|\nabla\chi_i^n\|_0\,\|\chi_i^n\|_0.\nonumber%\label{eq-i4}
 \end{eqnarray}
 Using  the  inequality $mn\leq \dfrac{m^2}{\theta}+n^2\theta\,$\,\,\text{with}\, $\theta:=\frac{\xi_{i,0}\,\min\lbrace\,1,\beta_*\,\rbrace}{8}$ and from \eqref{veminterpol}, we get
 
 \begin{eqnarray}
 I_4&\leq&(C\,\|(\omega_1,\omega_2)\|_{\infty,s})^2\,\frac{8}{\xi_{i,0}\,\min\lbrace\,1,\beta_*\,\rbrace}\,\delta t\,h^{2s+2}\,\|u_i^n\|^2_{s+1}+\frac{\xi_{i,0}\,\min\lbrace\,1,\beta_*\,\rbrace}{8}\,\delta t\,\|\nabla \chi_i^n\|^2_0\nonumber\\
 &&\,+\omega_{max}^2\,\delta t\,C\,h^{2s}\,\|u_i^n\|^2_{s+1}+\delta t\,\|\chi_i^n\|^2_0\nonumber\\
 & &\,+\delta t\,\frac{\xi_{i,0}\,\min\lbrace\,1,\beta_*\,\rbrace}{8}\,\|\nabla\chi_i^n\|^2_0+\delta t\,\frac{8\,\omega_{max}^2}{\xi_{i,0}\,\min\lbrace\,1,\beta_*\,\rbrace}\,\|\chi_i^n\|^2_0.\label{eq-i4}
 \end{eqnarray}

 Adding and subtracting the terms $\delta t\,\big( u_i^n\sum\limits_{j=1}^m\mathsf{A}(i,j)(\,u^{n-1}_j+\mathcal{P}u^{n-1}_j\,),\,\chi_i^n\, \big)_{0}$, \,$\delta t\,\big( \mathcal{P}u_i^n\sum\limits_{j=1}^m\mathsf{A}(i,j)\mathcal{P}U^{n-1}_j,\,\chi_i^n+\mathcal{P}\chi_i^n\, \big)_{0}$ and $\delta t\,\big( \mathcal{P}u_i^n\sum\limits_{j=1}^m\mathsf{A}(i,j)\mathcal{P}u^{n-1}_j,\,\chi_i^n\, \big)_{0}$, we get
 
 \begin{eqnarray}
I_5&:=&\delta t\,\Big( u^n_i\big(\textstyle\sum\limits_{j=1}^m\mathsf{A}(i,j)u^n_j\big),\,\chi_i^n\, \Big)_{0}-\delta t\,\Big( \mathcal{P}U^n_i\big(\textstyle\sum\limits_{j=1}^m\mathsf{A}(i,j)\mathcal{P}U^{n-1}_j\big),\,\mathcal{P}\chi_i^n\, \Big)_{0}\nonumber\\
&=&\delta t\,\Big( u^n_i\textstyle\sum\limits_{j=1}^m\mathsf{A}(i,j)\big(u^n_j-u^{n-1}_j\big),\,\chi_i^n\, \Big)_{0}+\delta t\,\Big( \mathcal{P}\varphi^n_i\big(\textstyle\sum\limits_{j=1}^m\mathsf{A}(i,j)\mathcal{P}U^{n-1}_j\big),\,\mathcal{P}\chi_i^n\, \Big)_{0}\nonumber\\
& &\,+\delta t\,\Big( u^n_i\textstyle\sum\limits_{j=1}^m\mathsf{A}(i,j)\big(u^{n-1}_j-\mathcal{P}u_j^{n-1}\big),\,\chi_i^n\, \Big)_{0}+\delta t\,\Big(\,(u^n_i-\mathcal{P}u^n_i)\textstyle\sum\limits_{j=1}^m\mathsf{A}(i,j)\mathcal{P}u^{n-1}_j,\,\chi_i^n\, \Big)_{0}\nonumber\\
& &\,+\delta t\,\Big( \mathcal{P}u^n_i\textstyle\sum\limits_{j=1}^m\mathsf{A}(i,j)\mathcal{P}\varphi_j^{n-1},\,\chi_i^n\, \Big)_{0}+\delta t\,\Big(\mathcal{P}u^n_i\textstyle\sum\limits_{j=1}^m\mathsf{A}(i,j)\mathcal{P}U^{n-1}_j,\,\chi_i^n-\mathcal{P}\chi_i^n\, \Big)_{0}\nonumber
 \end{eqnarray}
 
 Using generalised H$\Ddot{o}$lder's inequality and property of $\mathcal{P}$, we obtain,
 \begin{eqnarray}
 I_5&\leq& \delta t\,\mathsf{A}_{max}\,\Big[\,\|u_i^n\|_{\infty} \textstyle\sum\limits_{j=1}^m\|u^n_j-u^{n-1}_j\|_0+\|\varphi_i^n\|_{0} \textstyle\sum\limits_{j=1}^m\|\mathcal{P}U^{n-1}_j\|_\infty+\|u_i^n\|_{\infty} \textstyle\sum\limits_{j=1}^m\|u^{n-1}_j-\mathcal{P}u^{n-1}_j\|_0\nonumber\\
 & &\,+\|u_i^n-\mathcal{P}u_i^n\|_0 \textstyle\sum\limits_{j=1}^m\|\mathcal{P}u^{n-1}_j\|_\infty+\|\mathcal{P}u_i^n\|_{\infty} \textstyle\sum\limits_{j=1}^m\|\mathcal{P}\varphi_j^{n-1}\|_0\,\Big]\,\|\chi_i^n\|_0\nonumber\\
 & &\,+\delta t\,\mathsf{A}_{max}\,\|\mathcal{P}u_i^n\|_{\infty}\,\textstyle\sum\limits_{j=1}^m\big(\,\mathcal{P}U^{n-1}_j,\,\chi_i^n-\mathcal{P}\chi_i^n\, \big)_{0}\nonumber
 \end{eqnarray}
 
 Using Taylor's theorem w.r.t $t$ ( neglecting $(\delta t)^2$ term ), the estimates in \eqref{projlinf},\, \eqref{pi0interpol}, \eqref{veminterpol}, \,\eqref{assump}, \eqref{assump1} and then Young's inequality, we have,
 
 \begin{eqnarray}
 I_5&\leq&  \delta t\,\mathsf{A}_{max}\,\Big[\,\mu\, \textstyle\sum\limits_{j=1}^mC\,\delta t\,\|(u_j)^{n}_t\|_0+ m\,(\mu+1)\,(\|\eta_i^n\|_0+\|\chi_i^n\|_0)+\mu\,\textstyle\sum\limits_{j=1}^mC\,h^{s+1}\,\|u^n_j\|_{s+1}\nonumber\\
 & &\,+m\,(\mu+1)\,C\,h^{s+1}\,\|u^n_i\|_{s+1}+\mu\,\textstyle\sum\limits_{j=1}^m(\|\eta^{n-1}_j\|_0+\|\chi^{n-1}_j\|_0)\,\Big]\,\|\chi_i^n\|_0+0\nonumber\\%.\label{eq-i5}\
 &\leq&  \delta t\,\mathsf{A}_{max}\,\Big[\,\mu\, \textstyle\sum\limits_{j=1}^mC\,\delta t\,\|(u_j)^{n}_t\|_0+ m\,(\mu+1)\,(C\,h^{s+1}\,\|u_i^n\|_{s+1}+\|\chi_i^n\|_0)+\mu\,\textstyle\sum\limits_{j=1}^mC\,h^{s+1}\,\|u^n_j\|_{s+1}\nonumber\\
 & &\,+m\,(\mu+1)\,C\,h^{s+1}\,\|u^n_i\|_{s+1}+\mu\,\textstyle\sum\limits_{j=1}^m(C\,h^{s+1}\,\|u_i^{n-1}\|_{s+1}+\|\chi^{n-1}_j\|_0)\,\Big]\,\|\chi_i^n\|_0\nonumber\\
&\leq&
 \delta t\,\mathsf{A}_{max}\,\Big[\,\mu\,C\,\delta t\,\boldsymbol{\|u^{n}_t\|}_0+ m\,(\mu+1)\,(C\,h^{s+1}\,\|u_i^n\|_{s+1}+\|\chi_i^n\|_0)+\mu\,C\,h^{s+1}\,\boldsymbol{\|u^n\|}_{s+1}\nonumber\\
 & &\,+m\,(\mu+1)\,C\,h^{s+1}\,\|u^n_i\|_{s+1}+\mu\,C\,h^{s+1}\,\boldsymbol{\|u^{n-1}\|}_{s+1}+\mu\boldsymbol{\|\chi^{n-1}\|}_0\,\Big]\,\|\chi_i^n\|_0\nonumber\\
  &\leq& \mathsf{A}_{max}\,\Big[\,\mu\,C\,(\delta t)^2\, \boldsymbol{\|u^{n}_t\|}^2_0+ m\,(\mu+1)\,\delta t\,C\,h^{2s+2}\,\|u_i^n\|_{s+1}^2+\mu\,\delta t\,C\,h^{2s+2}\,\boldsymbol{\|u^n\|}^2_{s+1}\nonumber\\
 & &\,+m\,(\mu+1)\,\delta t\,C\,h^{2s+2}\,\|u^n_i\|^2_{s+1}+\mu\,\delta t\,C\,h^{2s+2}\,\boldsymbol{\|u^{n-1}\|}^2_{s+1} +\mu\,\delta t\,\boldsymbol{\|\chi^{n-1}\|}^2_0\nonumber\\
 & &\,+\delta t \, \big(\,4\mu+3m(\mu+1)\, \big)\,\|\chi_i^n\|^2_0\,\Big].\quad\label{eq-i5}
 \end{eqnarray}
 
 Now adding and subtracting $\delta t\,\left(\,\mathsf{R}(i,i)(\mathcal{P}u_i^n+\mathcal{P}U_i^n),\,\chi_i^n\, \right)_{0}$,\, using Cauchy-Schwarz inequality and property of $\mathcal{P}$, we get,
 
 \begin{eqnarray}
 I_6&:=&\delta t\,\left(\,\mathsf{R}(i,i)u_i^n,\,\chi_i^n\, \right)_{0}    -\delta t\, \left(\,\mathsf{R}(i,i)\mathcal{P}U^n_i,\,\mathcal{P}\chi_i^n\, \right)_{0}\,\nonumber\\
 &=&\delta t\,\left(\,\mathsf{R}(i,i)\,(u_i^n-\mathcal{P}u_i^n),\,\chi_i^n\, \right)_{0}+\delta t\,\left(\,\mathsf{R}(i,i)\,\mathcal{P}\varphi_i^n),\,\chi_i^n\, \right)_{0}+\delta t\, \left(\mathsf{R}(i,i)\mathcal{P}U^n_i,\,\chi_i^n-\mathcal{P}\chi_i^n\, \right)_{0}\nonumber\\
 &\leq& \delta t\,\mathsf{R}_{max}\,\|u_i^n-\mathcal{P}u_i^n\|_0\,\|\chi_i^n\|_{0}+\delta t\,\mathsf{R}_{max}\,\|\varphi_i^n\|_0\,\|\chi_i^n\|_{0}+\delta t\,\mathsf{R}_{max} \left(\mathcal{P}U^n_i,\,\chi_i^n-\mathcal{P}\chi_i^n\, \right)_{0}\nonumber\\
 &\leq& \delta t\,\mathsf{R}_{max}\,\big[\,C\,h^{s+1}\,\|u_i^n\|_{s+1}+\|\eta_i^n\|_0+\|\chi_i^n\|_0\,\big]\,\|\chi_i^n\|_{0}\qquad\quad\text{( use \eqref{pi0interpol}, \eqref{splite} )}\nonumber\\
 &\leq& \delta t\,\mathsf{R}_{max}\,\big[\,2\,C\,h^{s+1}\,\|u_i^n\|_{s+1}+\|\chi_i^n\|_0\,\big]\,\|\chi_i^n\|_{0}\qquad\quad\text{( use \eqref{veminterpol})}\nonumber %\label{eq-i6}
 \end{eqnarray}
 
 Using Young's inequality, we obtain
 \begin{eqnarray}
 I_6 &\leq&\mathsf{R}_{max}\,\Big[\, \delta t\,C\,4\,h^{2s+2}\,\|u_i^n\|^2_{s+1}+2\,\delta t\,\|\chi_i^n\|^2_{0}\,\big].\qquad\label{eq-i6}
 \end{eqnarray}
  Adding and subtracting $\big( \textstyle\sum\limits_{\ell,j\neq i}\mathsf{Q}(i,\ell,j)\mathcal{P}U^{n-1}_\ell \mathcal{P}U^{n-1}_j,\,\chi_i^n\, \big)_{0}$, we get
 \begin{eqnarray}
 I_7&=&\delta t\,\Big( \big(\textstyle\sum\limits_{\ell,j\neq i}\mathsf{Q}(i,\ell,j)u^n_\ell u^n_j\big),\,\chi_i^n\, \Big)_{0}-\delta t\,\Big( \big(\textstyle\sum\limits_{\ell,j\neq i}\mathsf{Q}(i,\ell,j)\mathcal{P}U^{n-1}_\ell \mathcal{P}U^{n-1}_j\big),\,\mathcal{P}\chi_i^n\, \Big)_{0}\nonumber\\
 &=&\delta t\,\Big( \textstyle\sum\limits_{\ell,j\neq i}\mathsf{Q}(i,\ell,j)\,\big(u^n_\ell u^n_j-\mathcal{P}U^{n-1}_\ell \mathcal{P}U^{n-1}_j\big),\,\chi_i^n\, \Big)_{0}\nonumber\\
 & &\,+\delta t\,\Big( \big(\textstyle\sum\limits_{\ell,j\neq i}\mathsf{Q}(i,\ell,j)\mathcal{P}U^{n-1}_\ell \mathcal{P}U^{n-1}_j\big),\,\chi_i^n-\mathcal{P}\chi_i^n\, \Big)_{0}.\label{eqt-i7}
 \end{eqnarray}
 
 We estimate
 \begin{eqnarray}
 u^n_\ell u^n_j-\mathcal{P}U^{n-1}_\ell \mathcal{P}U^{n-1}_j
 &=&u^n_\ell( u^n_j-\mathcal{P}U^{n-1}_j)+ \mathcal{P}U^{n-1}_j(u^n_\ell-\mathcal{P}U^{n-1}_\ell)\nonumber\\
  &=&u^n_\ell\big[\,( u^n_j- u^{n-1}_j)+(I-\mathcal{P})u^{n-1}_j+\mathcal{P}\varphi^{n-1}_j\,\big]\nonumber\\
  &\,&\quad+\mathcal{P}U^{n-1}_j\big[\,(u^n_\ell-u^{n-1}_\ell)+(I-\mathcal{P})u^{n-1}_\ell+\mathcal{P}\varphi^{n-1}_\ell\,\big].\qquad\qquad\label{eqt}
 \end{eqnarray}
 
 Combining \eqref{eqt-i7} and \eqref{eqt}, and using generalised H$\Ddot{o}$lder's inequality  we obtain
 \begin{eqnarray}
I_7&\leq&  \delta t\,\mathsf{Q}_{max}\,\textstyle\sum\limits_{\ell\neq i}\|u_\ell^n\|_\infty\textstyle\sum\limits_{j\neq i}\big[\,\| u^n_j- u^{n-1}_j\|_0+\|(I-\mathcal{P})u^{n-1}_j\|_0+\|\mathcal{P}\varphi^{n-1}_j\|_0\,\big]\,\|\chi_i^n\|_0\nonumber\\
& &\,+ \delta t\,\mathsf{Q}_{max}\,\textstyle\sum\limits_{j\neq i}\|\mathcal{P}U_j^{n-1}\|_\infty\textstyle\sum\limits_{\ell\neq i}\big[\,\| u^n_\ell- u^{n-1}_\ell\|_0+ \|(I-\mathcal{P})u^{n-1}_\ell\|_0 +\|\mathcal{P}\varphi^{n-1}_\ell \|_0\,\big] \,\|\chi_i^n\|_0 \nonumber\\
 & &\,+ \delta t\,\mathsf{Q}_{max}\,\textstyle\sum\limits_{\ell\neq i}\|\mathcal{P}U^{n-1}_\ell\|_\infty\textstyle\sum\limits_{j\neq i}\big(\, \mathcal{P}U^{n-1}_j,\,\chi_i^n-\mathcal{P}\chi_i^n\, \big)_{0}.\nonumber
 \end{eqnarray}
 
 Using Taylor's theorem w.r.t $t$ ( neglecting $(\delta t)^2$ term ), the estimates in \eqref{projlinf},\, \eqref{pi0interpol}, \eqref{splite}, \eqref{veminterpol}, \,\eqref{assump}, \eqref{assump1} and Young's inequality, we have,
 
 \begin{eqnarray}
 I_7&\leq& \delta t\,\mathsf{Q}_{max}\,C\,m\,\Big\lbrace\,\mu\,\textstyle\sum\limits_{j\neq i}\big[\,\delta t\,\| (u_j)^{n}_t\|_0+2\,h^{s+1}\,\|u^{n-1}_j\|_{s+1}+\|\chi^{n-1}_j\|_0\,\big]\,\|\chi_i^n\|_0\nonumber\\
 & &\,+(\mu+1)\,\textstyle\sum\limits_{\ell\neq i}\big[\,\delta t\,\| (u_\ell)^{n}_t\|_0+2\,h^{s+1}\,\|u^{n-1}_\ell\|_{s+1}+\|\chi^{n-1}_\ell\|_0\,\big]\,\|\chi_i^n\|_0\,\Big\rbrace+0\nonumber\\
 &\leq& 2\, (\delta t) \,\mathsf{Q}_{max}\,C\,m\,(\mu+1)\,\big\lbrace\,\delta t\,\boldsymbol{\| u^{n}_t\|}_0 +2\,h^{s+1}\,\boldsymbol{\|u^{n-1}\|}_{s+1}+\boldsymbol{\|\chi^{n-1}\|}_0\,\big\rbrace\,\|\chi_i^n\|_0\nonumber\\
&\leq& 2 \,\mathsf{Q}_{max}\,C\,m\,(\mu+1)\,\big\lbrace\,(\delta t)^2\,\boldsymbol{\| u^{n}_t\|}^2_0 +\delta t\,4\,h^{2s+2}\,\boldsymbol{\|u^{n-1}\|}^2_{s+1}+\delta t\,\boldsymbol{\|\chi^{n-1}\|}^2_0+3\,\delta t\,\|\chi_i^n\|^2_0\,\big\rbrace.\nonumber\\
&&\qquad\qquad\qquad\label{eq-i7}
 \end{eqnarray}
 
 Adding and subtracting $\Big( \textstyle\sum\limits_{j\neq i}\mathsf{R}(i,j)\mathcal{P}U^{n-1}_j,\,\chi_i^n\, \Big)_{0}$, \,using Cauchy-Schwarz inequality and triangle inequality,\, we get 
 \begin{eqnarray}
I_8&:=& \delta t\,\Big( \big(\textstyle\sum\limits_{j\neq i}^m\mathsf{R}(i,j)u^n_j\big),\,\chi_i^n\, \Big)_{0}-\delta t\,\Big( \big(\textstyle\sum\limits_{j\neq i}\mathsf{R}(i,j)\mathcal{P}U^{n-1}_j\big),\,\mathcal{P}\chi_i^n\, \Big)_{0}\nonumber\\
&=& \delta t\,\Big( \textstyle\sum\limits_{j\neq i}^m\mathsf{R}(i,j)\big(u^n_j-\mathcal{P}U^{n-1}_j\big),\,\chi_i^n\, \Big)_{0}+\delta t\,\Big( \big(\textstyle\sum\limits_{j\neq i}\mathsf{R}(i,j)\mathcal{P}U^{n-1}_j\big),\,\chi_i^n-\mathcal{P}\chi_i^n\, \Big)_{0}\nonumber\\
&\leq& \delta t\,\mathsf{R}_{max}\, \textstyle\sum\limits_{j\neq i}^m\|u^n_j-\mathcal{P}U^{n-1}_j\|_0\,\|\chi_i^n\|_{0}+\delta t\,\mathsf{R}_{max}\,\textstyle\sum\limits_{j\neq i}\big(\,\mathcal{P}U^{n-1}_j,\,\chi_i^n-\mathcal{P}\chi_i^n\, \big)_{0}\nonumber\\
&\leq& \delta t\,\mathsf{R}_{max}\, \textstyle\sum\limits_{j\neq i}^m\big[\,\|u^n_j-u^{n-1}_j\|_0+\|u^{n-1}_j-\mathcal{P}u^{n-1}_j\|_0+\|\mathcal{P}\varphi^{n-1}_j\|_0\,\big]\,\|\chi_i^n\|_{0} + 0.\nonumber
 \end{eqnarray}
 Using Taylor's theorem w.r.t $t$ ( neglecting $(\delta t)^2$ term ),  \eqref{pi0interpol}, property of  $\mathcal{P}$, \eqref{splite}, \eqref{veminterpol} and Young's inequality, we obtain,
 \begin{eqnarray}
 I_8&\leq& \delta t\,\mathsf{R}_{max}\,C\, \textstyle\sum\limits_{j\neq i}^m\big[\,\delta t\,\|(u_j)_t^{n}\|_0+2\,h^{s+1}\,\|u^{n-1}_j\|_{s+1}+\|\chi^{n-1}_j\|_0\,\big]\,\|\chi_i^n\|_{0} + 0\nonumber\\
% &\leq& \delta t\,\mathsf{R}_{max}\,C\, \big[\,\delta t\,\boldsymbol{\|u_t^{n}\|}_0+2\,h^{s+1}\,\boldsymbol{\|u^{n-1}\|}_{s+1}+\boldsymbol{\|\chi^{n-1}\|}_0\,\big]\,\|\chi_i^n\|_{0}\nonumber\\
  &\leq& \mathsf{R}_{max}\,C\, \big[\,(\delta t)^2\,\boldsymbol{\|u_t^{n}\|}^2_0+\delta t\,4\,h^{2s+2}\,\boldsymbol{\|u^{n-1}\|}_{s+1}^2+\delta t\,\boldsymbol{\|\chi^{n-1}\|}^2_0+3\,\delta t\,\|\chi_i^n\|^2_{0}\,\big].\qquad\qquad\label{eq-i8}
 \end{eqnarray}
 
 \noindent Note that $I_9:=(\delta t)\,\left(\,f^n_i,\,(I-\mathcal{P})\chi_i^n\,\right)_0 = (\delta t)\, \left(\,(I-\mathcal{P})f^n_i,\,\chi_i^n\,\right)_0$.\, Then using Cauchy-Schwarz inequality and Young's inequality, we have,
 \begin{eqnarray}
 I_9&\leq& \delta t\,\|(I-\mathcal{P})f^n_i\|_0\,\|\chi_i^n\|_0\,\leq\, \delta t\,C\,h^{s+1}\,\|f^n_i\|_{s+1}\,\|\chi_i^n\|_0\,\leq\,\delta t\,C\,h^{2s+2}\,\|f^n_i\|^2_{s+1}+\delta t\,\|\chi_i^n\|^2_0.\qquad\qquad\label{eq-i9}
 \end{eqnarray}
 
 Substituting the estimates \eqref{x3}, \eqref{x4}, \eqref{eq-i1}-\eqref{eq-i6} and \eqref{eq-i7}-\eqref{eq-i8} into  the equation \eqref{mx1}, and simplifying , we obtain
 \begin{eqnarray}
 &&\frac{\min\lbrace1,\alpha_*\rbrace}{4}\boldsymbol{\|\chi^n\|}_0^2+\frac{\xi_{min}\,\min\lbrace\,1,\beta_*\,\rbrace}{2}(\delta t)\,\boldsymbol{\|\nabla \chi^n\|}_0^2\nonumber\\
&&\leq \, \frac{2(1+\alpha^*)^2}{\min\lbrace1,\alpha_*\rbrace}\,\boldsymbol{\|\chi^{n-1}\|}_0^2+\mathcal{C}_1\,\delta t\,\boldsymbol{\|\chi^n\|}_0^2+\mathcal{C}_2\,\delta t\,\boldsymbol{\|\chi^{n-1}\|}_0^2+\mathcal{C}_3\,\delta t\,h^{2s}\,\boldsymbol{\|u^{n}\|}_{s+1}^2\nonumber\\
&&\quad\, \,+\mathcal{C}_4\,\delta t\,h^{2s}\,\boldsymbol{\|u^{n-1}\|}_{s+1}^2 +\mathcal{C}_5\,(\delta t)^2\,\boldsymbol{\|u_t^n\|}_0^2+(\delta t)^2\,\boldsymbol{\|u_{tt}\|}_0^2+\delta t\,h^{2s}\,\boldsymbol{\|f^n\|}_{s+1}^2,\qquad\qquad\label{first-ineq}
 \end{eqnarray}
 where,
 
 \begin{eqnarray}
 \mathcal{C}_1&:=&\frac{8\,\omega_{max}^2}{\xi_{min}\,\min\lbrace\,1,\beta_*\,\rbrace}+\mathsf{A}_{max}\big(\,4\mu+3m(\mu+1)\, \big)+\mathsf{R}_{max}(3+2C)+\mathsf{Q}_{max}Cm(\mu+1),\nonumber\\
 \mathcal{C}_2&:=&m\mu\mathsf{A}_{max}+C\mathsf{Q}_{max}(\mu+1)m^2+m\mathsf{R}_{max},\nonumber\\
 \mathcal{C}_3&:=&\frac{C16(\alpha^*+1)}{\min\lbrace1,\alpha_*\rbrace}+C\omega^2_{max}+C\mathsf{A}_{max}(2m+3m\mu)+C4\mathsf{R}_{max}\nonumber\\
 &&\,+\frac{8C}{\xi_{min}\,\min\lbrace\,1,\beta_*\,\rbrace}\Big[\,(1+\beta_*)+\max\limits_{i}\lbrace \big(\max\limits_E\Bar{\xi}_i^E+\|\frac{\partial^s\xi_i}{\partial x^s}\|_\infty+\beta^*\big)^2 \rbrace+C\,\|(\omega_1,\omega_2)\|_{\infty,s}^2\,\Big],\nonumber\\
 \mathcal{C}_4&:=&C\frac{16(\alpha^*+1)}{\min\lbrace1,\alpha_*\rbrace}+Cm\mu\mathsf{A}_{max}+C\mathsf{Q}_{max}m^2(\mu+1)+4Cm\mathsf{Q}_{max}\nonumber\\
 \mathcal{C}_4&:=&\frac{8(\alpha^*+1)}{\min\lbrace1,\alpha_*\rbrace}+C\mathsf{A}_{max}m\mu+C\mathsf{Q}_{max}m^2(\mu+1)+Cm\mathsf{R}_{max}.\nonumber
 \end{eqnarray}
 Summing from $n=1$ to $k$, and noting $\boldsymbol{\|\chi^0\|}_{0}=0$ ( choose $U^0=\mathcal{U}^0$ ) we get,
 \begin{eqnarray}
 &&\frac{\min\lbrace1,\alpha_*\rbrace}{4}\boldsymbol{\|\chi^k\|}_0^2+\frac{\xi_{min}\,\min\lbrace\,1,\beta_*\,\rbrace}{2}\textstyle \sum\limits_{n=1}^k(\delta t)\,\boldsymbol{\|\nabla \chi^n\|}_0^2\nonumber\\
&&\leq \, C\,\frac{2(1+\alpha^*)^2}{\min\lbrace1,\alpha_*\rbrace}\,\boldsymbol{\|\chi^{0}\|}_0^2+C\,(\mathcal{C}_1+\mathcal{C}_2)\,\textstyle\sum\limits_{n=1}^k\delta t\,\boldsymbol{\|\chi^n\|}_0^2+(\mathcal{C}_3+\mathcal{C}_4)\,h^{2s}\,\sum\limits_{n=1}^k\delta t\,\boldsymbol{\|u^{n}\|}_{s+1}^2\nonumber\\
&&\quad\, \,+(\delta t)^2\,\textstyle\sum\limits_{n=1}^k(\,\mathcal{C}_5\,\boldsymbol{\|u_t^n\|}_0^2+\boldsymbol{\|u_{tt}\|}_0^2\,)+h^{2s}\,\sum\limits_{n=1}^k\delta t\,\boldsymbol{\|f^n\|}_{s+1}^2\nonumber\\
&&\leq \widetilde{C}_1\,\textstyle\sum\limits_{n=1}^k\delta t\,\boldsymbol{\|\chi^n\|}_0^2+\widetilde{C}_2\,h^{2s}\,\big(\,\boldsymbol{[u]}_{0,s+1}^2+\boldsymbol{[f]}_{0,s+1}\,\big)+\widetilde{C}_3\,(\delta t)^2\,(\,\boldsymbol{[u_t]}_{\infty,0}^2+\boldsymbol{[u_{tt}]}_{\infty,0}^2\,),\qquad\qquad\label{final_eqn}
\end{eqnarray}
where $\widetilde{C}_i's$ are dependent on $\mathcal{C}_1-\mathcal{C}_5$ and are independent of $h$ and $ \delta t.$ 
From \eqref{final_eqn}, note that,
\begin{eqnarray}
&&\frac{\min\lbrace1,\alpha_*\rbrace}{4}\boldsymbol{\|\chi^k\|}_0^2\leq \widetilde{C}_1\,\textstyle\sum\limits_{n=1}^k\delta t\,\boldsymbol{\|\chi^n\|}_0^2+\widetilde{C}_2\,h^{2s}\,\big(\,\boldsymbol{[u]}_{0,s+1}^2+\boldsymbol{[f]}_{0,s+1}\,\big)+\widetilde{C}_3\,(\delta t)^2\,(\,\boldsymbol{[u_t]}_{\infty,0}^2+\boldsymbol{[u_{tt}]}_{\infty,0}^2\,).\nonumber
\end{eqnarray}

Thus by discrete Gronwall's lemma, we obtain,
\begin{eqnarray}
\boldsymbol{\|\chi^n\|}_0^2\leq C\,\Big[\,\widetilde{C}_2\,h^{2s}\,\big(\,\boldsymbol{[u]}_{0,s+1}^2+\boldsymbol{[f]}_{0,s+1}\,\big)+\widetilde{C}_3\,(\delta t)^2\,(\,\boldsymbol{[u_t]}_{\infty,0}^2+\boldsymbol{[u_{tt}]}_{\infty,0}^2\,)\,\Big].\qquad\quad\label{normineqchi0}
\end{eqnarray}

Now substituting \eqref{normineqchi0} into \eqref{final_eqn}, we get
\begin{eqnarray}
\textstyle\sum\limits_{n=1}^k(\delta t)\,\boldsymbol{\|\nabla \chi^n\|}_0^2\leq C\,\Big[\,\widetilde{C}_2\,h^{2s}\,\big(\,\boldsymbol{[u]}_{0,s+1}^2+\boldsymbol{[f]}_{0,s+1}\,\big)+\widetilde{C}_3\,(\delta t)^2\,(\,\boldsymbol{[u_t]}_{\infty,0}^2+\boldsymbol{[u_{tt}]}_{\infty,0}^2\,)\,\Big].\qquad\qquad\label{normineqgrdchi}
\end{eqnarray}

%Using \eqref{splite}, \eqref{normineqchi0} and \eqref{veminterpol}, we obtain
%\begin{eqnarray}
%\boldsymbol{[u^n-U^n]}_{\infty,0}^2\leq C\,\big( h^{2s}+(\delta t)^2\,\big).\nonumber
%\end{eqnarray}

Using \eqref{splite}, \eqref{normineqgrdchi} and \eqref{veminterpol}, we obtain
\begin{eqnarray}
\boldsymbol{[u^n-U^n]}_{0,1}^2\leq C\,\big( h^{2s}+(\delta t)^2\,\big).\nonumber
\end{eqnarray}
\end{proof}
Next we will prove that the induction assumption $\boldsymbol{\|U^n\|}_\infty\leq (\mu+1)$ holds, using the following proposition ( see sec.5.6 in \cite{vemmaxnorm} ) :
\begin{Prop}
Let $w\in V_h^p$ and using the assumptions on the elements $E\in\mathcal{T}_h$, we have the norm inequality,
\begin{eqnarray}
\|w\|_{\infty,E}\,\leq \,h^{-1}\,\|w\|_{0,E}. \label{infineqo}
\end{eqnarray}
\end{Prop}

\begin{Lemma}
Let us consider the discrete solution $\boldsymbol{U^n}$ of \eqref{vemform} at the $n^{th}$ time step along with estimate \eqref{normineqchi0}. Then
\begin{eqnarray}
\boldsymbol{\|U^n\|}_\infty\leq (\mu+1).\label{induc_arg}
\end{eqnarray}
\end{Lemma}
\begin{proof}
From \eqref{splite}, \eqref{assump}, \eqref{infineqo}, \eqref{veminterpol} and \eqref{normineqchi0}, \, we have
\begin{eqnarray}
\boldsymbol{\|U^n\|}_\infty&\leq&\boldsymbol{\|U^n-u^n\|}_\infty+\boldsymbol{\|u^n\|}_\infty\leq \boldsymbol{\|\eta^n\|}_\infty+h^{-1}\,\boldsymbol{\|\chi^n\|}_0 +\mu\nonumber\\
&\leq&C\,(\,h^s+h^{s-1}+h^{-1}\,\delta t\,)\,\boldsymbol{ [u]}_{\infty,s+1}\,+\,\mu.\nonumber
\end{eqnarray}
For sufficiently small $h$ and $\delta t$, we obtain \eqref{induc_arg}.
\end{proof}

\section{Implementation details}
We outline two techniques to solve the system of equations obtained from the discrete formulation \eqref{vemformlinear}.

\noindent\textbf{1. \,Iteration method}

At each time step, to enhance the accuracy of our numerical solution, we solve the discrete linear system \eqref{vemformlinear} by using an iteration method. In $Algorithm\,1$, the iteration procedure is detailed in steps (2.2)-(2.3), and the advantages in its implementation is described in \cite{fembiyue}.
%To circumvent the difficulty of solving the nonlinear system of equations, we linearised the nonlinear terms in \eqref{vemform} at each time step and obtained the linear system \eqref{vemformlinear}. In the numerical execution of \eqref{vemformlinear}, we introduce an iteration procedure at each time level as outlined in $Algorithm \,1$

\noindent \textbf{\textit{Algorithm} 1 :}
\begin{itemize} 
    \item[(1)] Let $\boldsymbol{U^{0}} := \boldsymbol{u^{0}}$ on the  virtual element space $V_h^p$  and fix $tol>0$.
    \item[(2)] START FOR\,\, n=1,2,...,N \,\,\,DO
    \item[(2.1)] Set $\boldsymbol{\mathcal{N}^{\,0}} := \boldsymbol{U^{n-1}}$ and let $r=0$.
    \item[(2.2)] Set $r=r+1$. Find $\boldsymbol{\mathcal{N}^{\,r}}\in \boldsymbol{\mathcal{V}_h^p}$\,\,\, satisfing \,\,$ \forall\,\boldsymbol{V}\in \boldsymbol{\mathcal{V}_h^p}$  : 
\begin{eqnarray}
&&\textstyle\sum\limits_{i=1}^m\Big\lbrace\,m_h(\mathcal{N}_i^{\,r},V_i)-m_h(U_i^{n-1},V_i)+(\delta t)\,a_h(\mathcal{N}_i^{\,r},V_i)+(\delta t)\,l_{1,h}(\mathcal{N}_i^{\,r},V_i)+(\delta t)\,l_{2,h}(\mathcal{N}_i^{\,r},V_i) \nonumber\\
&&\quad+(\delta t)
\,l_{3,h}(\mathcal{N}_i^{\,r},V_i)\Big\rbrace \, =\, \textstyle (\delta t)\sum\limits_{i=1}^m\Big\lbrace-
l_{4,h}(\mathcal{N}_i^{\,r},V_i)- l_{5,h}(\mathcal{N}_i^{\,r},V_i)+ (\,f^n_i,\,\mathcal{P}V_i\,)_0\,\Big\rbrace.\nonumber
\end{eqnarray}
\item[(2.3)] Repeat steps (2.2) UNTIL $ \,\textit{norm }( \boldsymbol{\mathcal{N}^{\,r}-\mathcal{N}^{\,r-1}})\,<\,tol$.
\item[(2.4)] Set $ \boldsymbol{U^n}:=\boldsymbol{\mathcal{N}^{\,r}}$.\quad END FOR \hfill $\diamond$
\end{itemize}

\noindent\textbf{2.\,Two-grid method }

 Usually two-grid method involves two meshes $\mathcal{T}_H,\,\mathcal{T}_h$ of different mesh diameters $H$, $h$ ( with $H<h$ ) along with the corresponding VEM spaces $V_H^p$ and $V_h^p$, known as coarse space and fine space, respectively. To our knowledge, this is the first instance where we use a two-grid concept to solve the discrete scheme obtained by VEM discretization of a system of equations.  
 Now we present our two-grid method to solve the virtual element formulation  \eqref{vemformlinear}.
%In our two-grid algorithm, for each time level, at the coarse mesh $\mathcal{T}_H$ using $Algorithm 1$ with relatively large $tol$,  we obtain a solution $U_H\in V_H^p$. Interpolate the coarse solution $U_H$ from $V_H^p$ to fine space $V_h^p$, use it as initial guess and perform few fixed number iterations in $Algorithm 1$ at the fine level solve. Now we present the two-grid scheme.

\noindent\textbf{\textit{Algorithm} 2 :}
\begin{itemize} 
    \item[(a)] Fix $ctol>0$ and an integer $fiter>0$.
    \item[(b)] START FOR\,\, n=1,2,...,N \,\,\,DO
    \item[(b.1)] On coarse mesh $\mathcal{T}_H$, \,\,find $\boldsymbol{\widetilde{U}}^n\in \boldsymbol{\mathcal{V}_H^p}$ solving \eqref{vemformlinear} using the steps (2.1)-(2.4) in $Algorithm 1$ by fixing $tol:=ctol$. 
    \item[(b.2)] Interpolate $\boldsymbol{\widetilde{U}}^n\in \boldsymbol{\mathcal{V}_H^p}$ to $ \boldsymbol{\mathcal{V}_h^p}$ and denote as $\boldsymbol{{\mathcal{I} U}}^n$.
    \item[(b.3)] Set $\boldsymbol{\mathcal{N}^{\,0}}:=\boldsymbol{{\mathcal{I} U}}^n$ and $r=0$. On fine mesh $\mathcal{T}_h$, execute $fiter$ times the step (2.2) of $Algorithm\,1$,  and obtain the solution $ \boldsymbol{U^n}\in \boldsymbol{\mathcal{V}_h^p}$.\quad END FOR \hfill $\diamond$
\end{itemize}
In $Algorithm\,2$, we can choose $ctol$ relatively larger,say $10^{-3}$ and a small $fiter$, say 1,2 or 3. In the section \ref{sec:numex}, we shall see that for the optimal performance of  our two-grid technique, the choice of $H$ wrt $h$ and the number of iterations $fiter$ in fine space, are problem dependent. 

\section{Numerical examples}
\label{sec:numex}
\begin{figure}[H]
	\subfloat[Distorted squares.] {{\includegraphics[width=5cm]{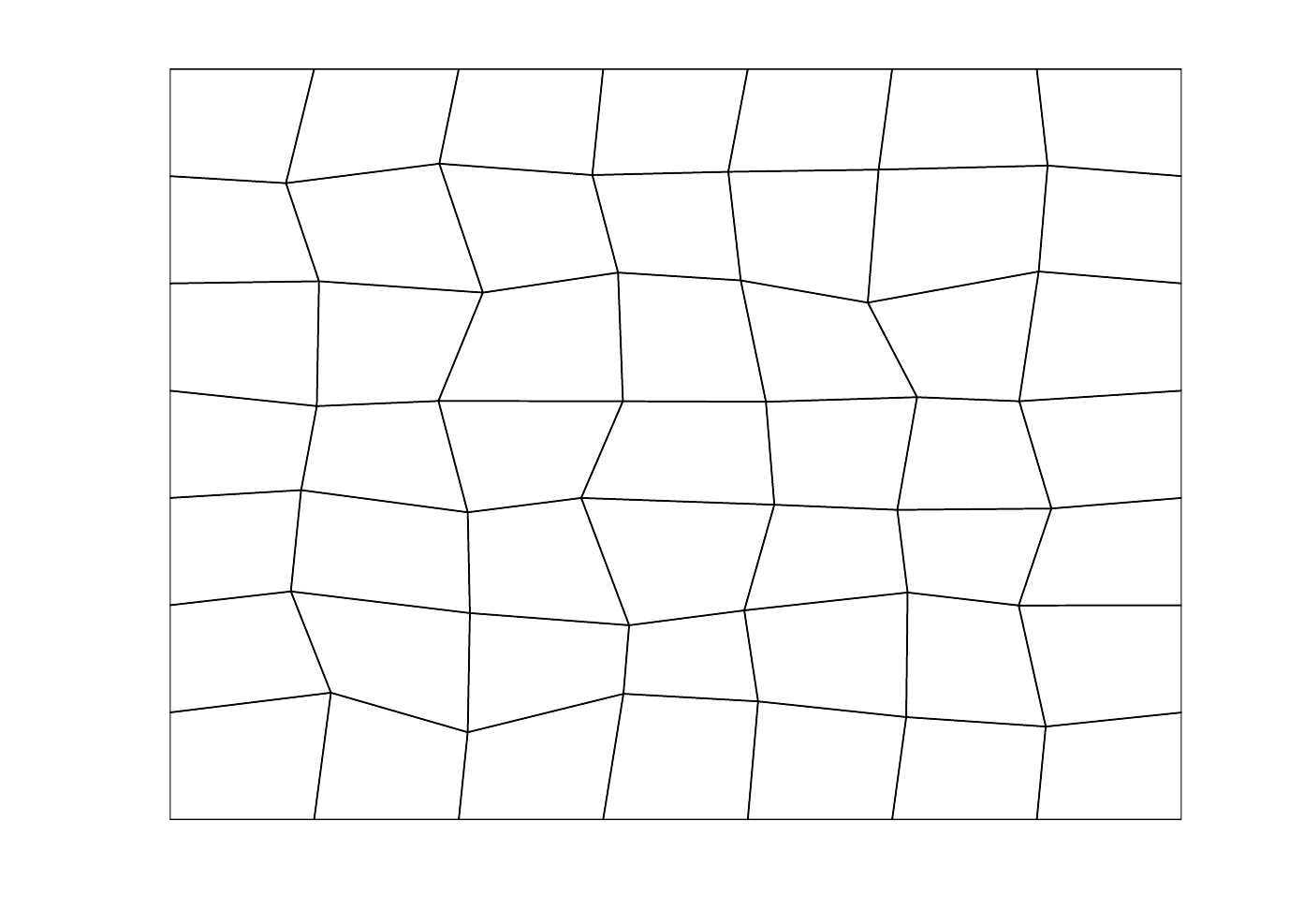} }}%
	\subfloat[Nonconvex mesh.] 
	{{\includegraphics[width=5cm]{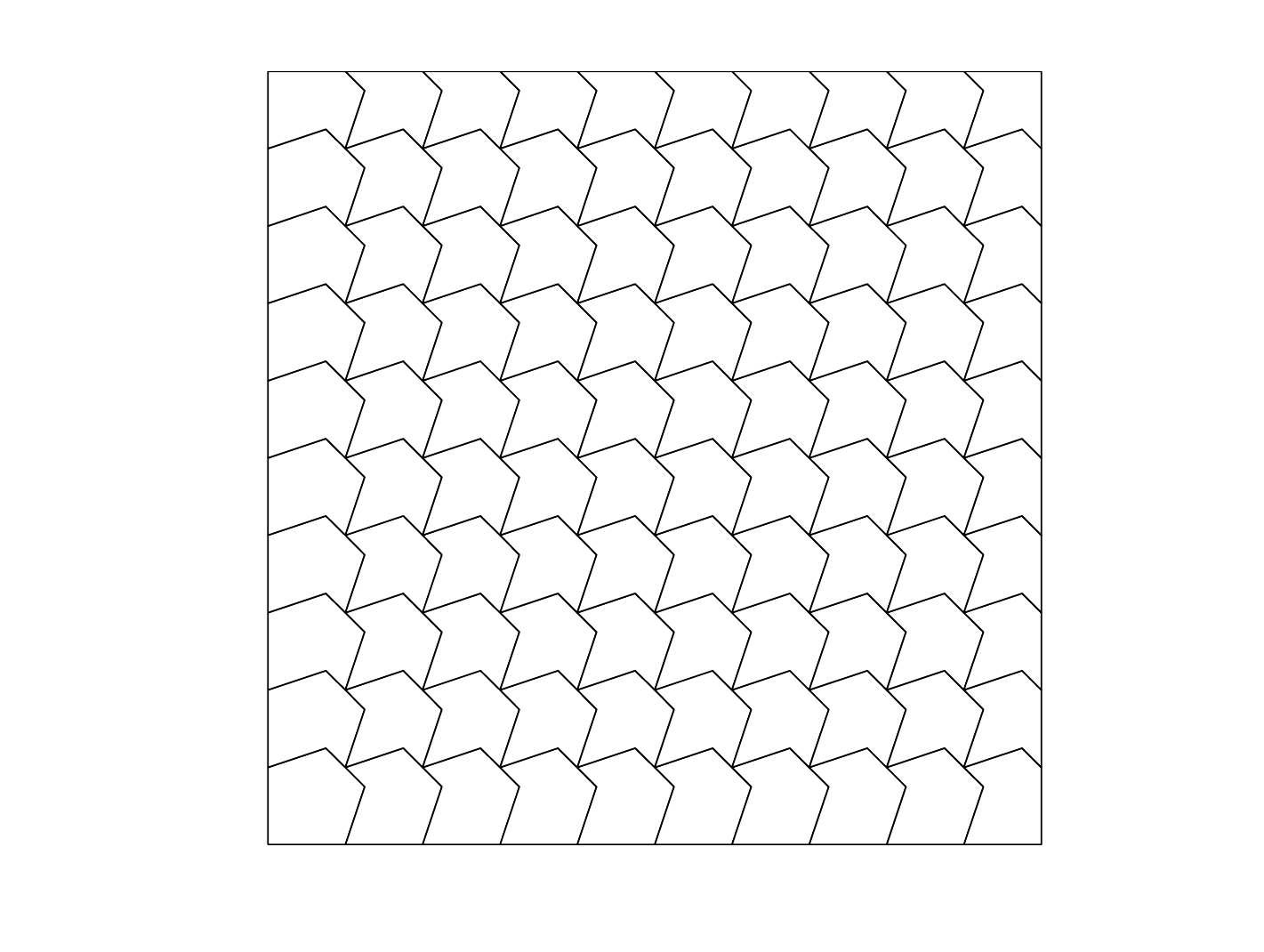} }}
	\subfloat[Regular Voronoi mesh.]
{{\includegraphics[width=4cm]{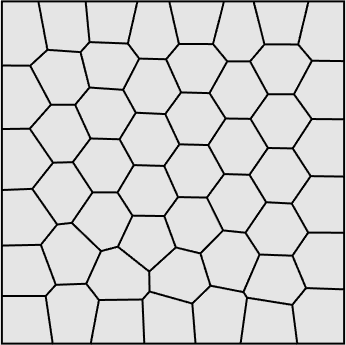} }}
	\caption{Samples of meshes with diameter $h=1/5$.}
	\label{ex2sampmesh} 
\end{figure}

In this section we consider two problems whose exact solutions are known.  For the numerical experiment, we use three types of meshes, namely, distorted  square mesh, nonconvex mesh and regular voronoi mesh. The choice of meshes include convex and concave elements.  A sample of each mesh is shown in figure \ref{ex2sampmesh}. 

To obtain the error, we evaluate norm of the difference of exact  and numerical solution, at the final time step $T$. Let $u$ and $u_h$ denote the exact and discrete solution,respectively, at time  $T$.  Then the $L^2$ norm and $H^1$ semi-norm of the error at $T$, denoted by $e_{h,0}$ and $e_{h,1}$ respectively, are calculated using the expressions,
\begin{eqnarray}
e_{h,0}^2 &=&\textstyle \sum\limits_{E\in\mathcal{T}_h} \|u - \Pi^0_pu_h \|_E^2 , \quad e_{h,1}^2 = \sum\limits_{E\in\mathcal{T}_h} \| \nabla (u - \Pi^{\nabla}_pu_h ) \|_E^2. \nonumber \nonumber 
\end{eqnarray}

The first example consists of a system of two equations, for which we show the rate of convergence of $H^1$ semi-norm with respect to $h$, for the VEM orders 1,2 and 3. The second example is a system of four equations and in this we show rate of convergence wrt $t$ for second-order VEM. On both the examples, , for VEM order k=2,  we compare the performances of Iteration method ( $Algorithm 1$ ) and two-grid method ( $Algorithm 2$ ),  based on accuracy in terms of error values $e_{h,0},\,e_{h,1}$ and on CPU time. 
 
\subsection{Example 1}
Let $m=2$, \,$\Omega=[0,1]\times[0,1]$ and $T=1$. \,The choice for parameters in \eqref{modelprob} are  $\,\xi_1=1,\,\xi_2=2,\,(\omega_1,\omega_2):=(1,2),\,\mathsf{Q}=[0],\,\mathsf{A}=[1,1.5;1.1,2],\,\mathsf{R}=[-1,0;2,0]$ and the functions $f_1$, $f_2$ are such that the exact solution $u_1$, $u_2$ are defined as :
\begin{eqnarray}
u_1(x,y,t):=e^t\,x\,y\,(x-1)^2\,(y-1)^2\quad\text{and}\quad
u_2(x,y,t):=e^{-t}\,x\,y\,(x-1)\,(y-1).\nonumber
\end{eqnarray}
Using $Algorithm\,1$, we solve the system of equations arising from the scheme \eqref{vemformlinear} for VEM orders k = 1,2 and 3, over the distorted square mesh. We choose $\delta t\approx \mathcal{O}(h^{k+1})$ and tolerance value $ tol:=10^{-7}$. The error $e_{h,1}$ for the solutions $u_1,\,u_2$ along  with rate of convergence ($roc$) are tabulated in table \ref{tab:e1h1err}. We see optimal rate wrt $H^1$ seminorm which is in accordance with our theoretical estimate.
\input{ex1mesh}

Now for VEM order k=2, we present a comparison of iteration method and two-grid method. For this sake, we fix $\delta t = 10^{-3}$. In $Algorithm\,1$, \,$tol$ is set to $10^{-6}$. In the two-grid scheme, we choose $H=2\,h$, set the coarse mesh solve tolerance $ctol$ as $10^{-3}$ and  perform one fine grid iteration, i.e. set $fiter=1$. 
The results obtained for distorted square mesh and nonconvex mesh can be found in table \ref{tab:e1h1err1-1} and table \ref{tab:e1h1err1-2}, respectively.
\input{compareeg11}
 We note that on both the meshes, the two-grid method performs efficiently by consuming considerably less CPU time than the standard scheme, without compromising on the accuracy. 

\subsection{Example 2}
Consider the example given in \cite{fembiyue}. We take $T=1$, $\Omega =[0,1] \times [0,1]$ and $m=4$. Let
$\xi_1=1$, $\xi_2=2$, $\xi_3=1.5$, $\xi_4=3$, $(\omega_1,\omega_2):=(1,2)$, 
$\mathsf{Q}=[0]$,

$\mathsf{A}=[1,1.5,0.5,1;1.1,2,0.7,1.2;1,1.5,3,0.4;1.5,0.5,1.4,2]$ and 

$\mathsf{R}=[-1,0,0,0;0,2,0,0;0,-1,0;0,0,0,-1]$.
The choice of coefficient for \eqref{modelprob} mimics the Lotka-Volterra competition type of reactions with the functions $f_i's,\,i=1,2,3,4$ defined as :
\begin{eqnarray}
&& u_1(x,y,t):= 10\, \text{sin}(\pi x)\,\text{sin}(\pi y)\,\text{cos}(t), \qquad
u_2(x,y,t):=5\,\text{sin} \Big(\frac{\pi}{2}x\Big)\, \text{sin}(\pi y)\,(1+\text{sin}(t)),  \nonumber\\
&& u_3(x,y,t):= 7\,\text{sin}(\pi x)\,\text{cos} \Big(\frac{\pi}{2}y\Big)\,\text{cos}(t), \qquad 
u_4(x,y,t):=12\,\text{sin}(\pi x)\,\text{sin}(\pi y)\,(x+y)\,(1+t).  \nonumber
\end{eqnarray}

To this example, we first show the time rate of convergence for second order VEM over regular Voronoi mesh. Thus, we fix $h=\frac{1}{16}$ and vary $\delta t = \frac{1}{5},\,\frac{1}{10},\,\frac{1}{20},\,\frac{1}{40}.$ The errors $e_{h,0}$ in table \ref{tab:e2roc-1}, were calculated  using the numerical solutions obtained from $Algorithm\,1$ with $tol=10^{-6}$.  The optimal rate of one for Euler discretization was obtained as seen in table \ref{tab:e2roc-1}.
\input{timeroc}

Next, let us fix $\delta t = 10^{-3}$. In the iteration method, set $tol$ to $10^{-6}$. For the two-grid method, we choose $H=\frac{1}{4}$ for $h=\frac{1}{8},\,\frac{1}{16}$ and $H=\frac{1}{8}$ for $h=\frac{1}{32}$, set $ctol$ as $10^{-3}$ and  perform three fine grid iteration, i.e. set $fiter=3$. We compute the results for second order VEM on regular Voronoi mesh (tables \ref{tab:e2h1err1-1}-\ref{tab:e2h1err1-11}) and  nonconvex mesh (tables \ref{tab:e2h1err1-2}-\ref{tab:e2h1err1-22}). 

\input{ex2comp}
\input{ex2comp1}
A glance of the tables \ref{tab:e2h1err1-1}-\ref{tab:e2h1err1-22} shows that for both the meshes, the accuracy of $e_{h,0},\,e_{h,1}$ on standard scheme and two-grid schemes remains the same. Table \ref{tab:e2h1err1-11t} clearly shows the computational efficiency of two-grid method. In particular, the effectiveness of two-grid method is more for smaller $h$ ( ie. for large degrees of freedom ).

\section{Conclusion}
In this work we have studied the virtual element method for a system of nonlinear advection-diffusion-reaction equation. To circumvent the difficulty of solving a nonlinear system, we have modified it to a linear discrete formulation whose solution also satisfies the model problem. For solving the arising system of equations, we propose iteration method and two-grid method. From the two examples, we clearly see the computational efficiency of two-grid method over iteration method. We also note that the choice of parameters $H$ and $fiter$ in $Algorithm\,2$ are problem dependent. 
\bibliographystyle{unsrt}
\bibliography{VEM_SUPG.bib}

\end{document}

%% file: ex1mesh.tex
\begin{table}[h]
\small
\centering % centering table
\begin{tabular}{|c| c c | c c | c c|} % creating 10 columns
\hline % inserting double-line
\multicolumn{7}{|c|}{$u_1$}\\
\hline
&\multicolumn{2}{|c|}{k=1} & \multicolumn{2}{|c|}{k=2} & \multicolumn{2}{|c|}{k=3} \\
\hline % inserts single-line
$h$  & $e_{h,1}$ & $roc$  & $e_{h,1}$ & $roc$ & $e_{h,1}$ & $roc$   \\ \hline
$ \frac{1}{4}$ & $3.73e^{-2}$  & $*$  & $5.28e^{-3}$ & $*$ & $5.69e^{-4}$ & $*$  \\ \hline
$\frac{1}{8}$ & $1.92e^{-2}$  & $0.96$  & $1.41e^{-3}$ & $1.90$ & $7.54e^{-5}$ & $2.91$ \\ \hline
$\frac{1}{16}$ & $9.69e^{-3}$  & $0.99$  & $3.61e^{-4}$ & $1.96$ & $9.96e^{-6}$ & $2.92$ \\ \hline
$\frac{1}{32}$ & $4.87e^{-3}$  & $0.99$  & $9.11e^{-5}$ & $1.99$ & $1.28e^{-6}$ & $2.96$  \\ \hline
\multicolumn{7}{|c|}{$u_2$}\\
\hline
&\multicolumn{2}{|c|}{k=1} & \multicolumn{2}{|c|}{k=2} & \multicolumn{2}{|c|}{k=3} \\
\hline % inserts single-line
$h$  & $e_{h,1}$ & $roc$  & $e_{h,1}$ & $roc$ & $e_{h,1}$ & $roc$ \\ \hline
$ \frac{1}{4}$& $1.08e^{-2}$ & $*$  & $1.13e^{-3}$ & $*$ & $7.06e^{-5}$ & $*$ \\ \hline
$\frac{1}{8}$ & $5.45e^{-3}$ & $0.98$  & $2.91e^{-4}$ & $1.96$ & $9.01e^{-6}$ & $2.97$ \\ \hline
$\frac{1}{16}$ & $2.75e^{-3}$ & $0.98$  & $7.34e^{-5}$ & $1.99$ & $1.26e^{-6}$ & $2.84$ \\ \hline
$\frac{1}{32}$  & $1.38e^{-3}$ & $0.99$  & $1.85e^{-5}$ & $1.99$ & $1.53e^{-7}$ & $3.04$ \\ \hline
\end{tabular}
\caption{Error $e_{h,1}$ for $u_1$ and $u_2$, along with $roc$ for distorted square mesh.} % title name of the table
\label{tab:e1h1err}
\end{table}

%% file: compareeg11.tex
\begin{table}[h]
\small
\centering % centering table
\begin{tabular}{|c| c| c c|  c c|} % creating 10 columns
\hline % inserting double-line
\multicolumn{6}{|c|}{ k=2. Iteration method.}\\
\hline
& Time &\multicolumn{2}{|c|}{$u_1$} & \multicolumn{2}{|c|}{$u_2$} \\\cline{3-6}
$h$&(sec)  &  $e_{h,0}$  & $e_{h,1}$  & $e_{h,0}$  & $e_{h,1}$   \\ \hline
$\frac{1}{8}$&73& $1.277061e^{-5}$  & $1.406609e^{-3}$  &  $2.565228e^{-6}$ & $2.906023e^{-4}$ \\ \hline
$\frac{1}{16}$&354&$1.744696e^{-6}$  & $3.611398e^{-4}$  &  $3.666301e^{-7}$ & $7.343147e^{-5}$ \\ \hline
$\frac{1}{32}$&2375& $6.207318e^{-7}$  & $9.115492e^{-5}$  & $1.707884e^{-7}$ & $1.850523e^{-5}$ \\ \hline
\multicolumn{6}{|c|}{ k=2. Two-grid method.}\\
\hline
& Time &\multicolumn{2}{|c|}{$u_1$} & \multicolumn{2}{|c|}{$u_2$} \\\cline{3-6}
$h$&(sec)  &  $e_{h,0}$  & $e_{h,1}$  & $e_{h,0}$  & $e_{h,1}$   \\ \hline
$\frac{1}{8}$&50& $1.281077e^{-5}$  & $1.406618e^{-3}$  &  $2.593940e^{-6}$ & $2.906076e^{-4}$ \\ \hline
$\frac{1}{16}$&227&$1.826598e^{-6}$  & $3.611482e^{-4}$  &  $3.956566e^{-7}$ & $7.343445e^{-5}$ \\ \hline
$\frac{1}{32}$&1576& $6.730950e^{-7}$  & $9.116259e^{-5}$  & $1.84544e^{-7}$ & $1.850782e^{-5}$ \\ \hline
\end{tabular}
\caption{CPU  time and error values $e_{h,0},\,e_{h,1}$ for distorted square mesh.} % title name of the table
\label{tab:e1h1err1-1}
\end{table}

\begin{table}[h]
\small
\centering % centering table
\begin{tabular}{|c| c| c c|  c c|} % creating 10 columns
\hline % inserting double-line
\multicolumn{6}{|c|}{ k=2. Iteration method.}\\
\hline
& Time &\multicolumn{2}{|c|}{$u_1$} & \multicolumn{2}{|c|}{$u_2$} \\\cline{3-6}
$h$&(sec)  &  $e_{h,0}$  & $e_{h,1}$  & $e_{h,0}$  & $e_{h,1}$   \\ \hline
$\frac{1}{8}$&139& $1.381993e^{-5}$  & $1.530203e^{-3}$  &  $2.534569e^{-6}$ & $2.964664e^{-4}$ \\ \hline
$\frac{1}{16}$&751&$1.760424e^{-6}$  & $3.795994e^{-4}$  &  $3.517047e^{-7}$ & $7.426381e^{-5}$ \\ \hline
$\frac{1}{32}$&10332& $6.176535e^{-7}$  & $9.425080e^{-5}$  & $1.697609e^{-7}$ & $1.858693e^{-5}$ \\ \hline
\multicolumn{6}{|c|}{ k=2. Two-grid method.}\\
\hline
& Time &\multicolumn{2}{|c|}{$u_1$} & \multicolumn{2}{|c|}{$u_2$} \\\cline{3-6}
$h$&(sec)  &  $e_{h,0}$  & $e_{h,1}$  & $e_{h,0}$  & $e_{h,1}$   \\ \hline
$\frac{1}{8}$&88& $1.381199e^{-5}$  & $1.530202e^{-3}$  &  $2.539033e^{-6}$ & $2.964676e^{-4}$ \\ \hline
$\frac{1}{16}$&509&$1.793708e^{-6}$  & $3.796026e^{-4}$  &  $3.639653e^{-7}$ & $7.426499e^{-5}$ \\ \hline
$\frac{1}{32}$&6021& $6.416986e^{-7}$  & $9.425410e^{-5}$  & $1.759072e^{-7}$ & $1.858805e^{-5}$ \\ \hline
\end{tabular}
\caption{ CPU  time and error values $e_{h,0},\,e_{h,1}$ for nonconvex mesh.} % title name of the table
\label{tab:e1h1err1-2}
\end{table}

%% file: timeroc.tex
\begin{table}[h]
\small
\centering % centering table
\begin{tabular}{|c|c|c||c|c||c|c||c|c|} % creating 10 columns
\hline
&\multicolumn{2}{|c||}{$u_1$} & \multicolumn{2}{|c||}{$u_2$}&\multicolumn{2}{|c||}{$u_3$} & \multicolumn{2}{|c|}{$u_4$} \\\cline{2-9}
$\delta t$  &  $e_{h,0}$  & $roc$&  $e_{h,0}$  & $roc$&  $e_{h,0}$  & $roc$&  $e_{h,0}$  & $roc$  \\ \hline
$\frac{1}{5}$&  $5.67e^{-3}$  & $*$&  $1.62e^{-3}$  & $*$&  $2.79e^{-3}$  & $*$&  $1.59e^{-3}$  & $*$  \\ \hline
$\frac{1}{10}$&  $2.71e^{-3}$  & $1.06$&  $8.47e^{-4}$  & $0.93$&  $1.33e^{-3}$  & $1.06$&  $8.06e^{-4}$  & $0.98$  \\ \hline
$\frac{1}{20}$&  $1.32e^{-3}$  & $1.04$&  $4.33e^{-4}$  & $0.97$&  $6.51e^{-4}$  & $1.03$&  $4.69e^{-4}$  & $0.78$  \\ \hline
$\frac{1}{40}$&  $6.55e^{-4}$  & $1.01$&  $2.22e^{-4}$  & $0.96$&  $3.22e^{-4}$  & $1.02$&  $3.44e^{-4}$  & $0.45$  \\ \hline
\end{tabular}
\caption{Rate of convergence wrt $\delta t$ for Voronoi mesh.} % title name of the table
\label{tab:e2roc-1}
\end{table}

%% file: ex2comp.tex
\begin{table}[H]
\small
\centering % centering table
\begin{tabular}{| c| c c||  c c|} % creating 10 columns
\hline % inserting double-line
\multicolumn{5}{|c|}{ Iteration method.}\\
\hline
& \multicolumn{2}{|c||}{$u_1$} & \multicolumn{2}{|c|}{$u_2$} \\\cline{2-5}
$h$&  $e_{h,0}$  & $e_{h,1}$  & $e_{h,0}$  & $e_{h,1}$   \\ \hline
$\frac{1}{8}$& $4.256581e^{-4}$  & $4.642544e^{-2}$  &  $3.562504e^{-4}$ & $3.825449e^{-2}$ \\ \hline
$\frac{1}{16}$&$5.554100e^{-5}$  & $1.091555e^{-2}$  &  $4.261672e^{-5}$ & $9.094519e^{-3}$ \\ \hline
$\frac{1}{32}$& $2.650741e^{-5}$  & $2.670695e^{-3}$  & $1.014168e^{-5}$ & $2.238159e^{-3}$ \\ \hline
&\multicolumn{2}{|c||}{$u_3$} & \multicolumn{2}{|c|}{$u_4$} \\\cline{2-5}
$h$  &  $e_{h,0}$  & $e_{h,1}$  & $e_{h,0}$  & $e_{h,1}$   \\ \hline
$\frac{1}{8}$& $1.497008e^{-4}$  & $1.586691e^{-2}$  &  $2.454550e^{-3}$ & $2.765068e^{-1}$ \\ \hline
$\frac{1}{16}$&$2.138141e^{-5}$  & $3.755053e^{-3}$  &  $2.928331e^{-4}$ & $6.665600e^{-2}$ \\ \hline
$\frac{1}{32}$& $1.287428e^{-5}$  & $9.125372e^{-4}$  & $3.569207e^{-5}$ & $1.602872e^{-2}$ \\ \hline
\end{tabular}
\caption{Error values $e_{h,0},\,e_{h,1}$ for Voronoi mesh.} % title name of the table
\label{tab:e2h1err1-1}
\end{table}

\begin{table}[H]
\small
\centering % centering table
\begin{tabular}{| c| c c||  c c|} % creating 10 columns
\hline % inserting double-line
\multicolumn{5}{|c|}{ Two-grid method.}\\
\hline
& \multicolumn{2}{|c||}{$u_1$} & \multicolumn{2}{|c|}{$u_2$} \\\cline{2-5}
$h$&  $e_{h,0}$  & $e_{h,1}$  & $e_{h,0}$  & $e_{h,1}$   \\ \hline
$\frac{1}{8}$& $4.249176e^{-4}$  & $4.642545e^{-2}$  &  $3.561006e^{-4}$ & $3.825466e^{-2}$ \\ \hline
$\frac{1}{16}$&$5.021513e^{-5}$  & $1.091543e^{-2}$  &  $4.472873e^{-5}$ & $9.094471e^{-3}$ \\ \hline
$\frac{1}{32}$& $2.244500e^{-5}$  & $2.670014e^{-3}$  & $6.332434e^{-6}$ & $2.238086e^{-3}$ \\ \hline
&\multicolumn{2}{|c||}{$u_3$} & \multicolumn{2}{|c|}{$u_4$} \\\cline{2-5}
$h$  &  $e_{h,0}$  & $e_{h,1}$  & $e_{h,0}$  & $e_{h,1}$   \\ \hline
$\frac{1}{8}$& $1.493283e^{-4}$  & $1.586688e^{-2}$  &  $2.454474e^{-3}$ & $2.765073e^{-1}$ \\ \hline
$\frac{1}{16}$&$1.876669e^{-5}$  & $3.755002e^{-3}$  &  $3.033806e^{-4}$ & $6.665607e^{-2}$ \\ \hline
$\frac{1}{32}$& $1.165073e^{-5}$  & $9.124825e^{-4}$  & $4.246686e^{-5}$ & $1.602876e^{-2}$ \\ \hline
\end{tabular}
\caption{Error values $e_{h,0},\,e_{h,1}$ for Voronoi mesh.} % title name of the table
\label{tab:e2h1err1-11}
\end{table}

%\begin{table}[h]
%\centering % centering table
%\begin{tabular}{| c|  c|   c|} % creating 10 columns
%\hline
%& Standard scheme & Two-grid scheme \\
%h& Time (sec)&Time(sec)\\ \hline
%$\frac{1}{8}$& 387  &315\\ \hline
%$\frac{1}{16}$&2842  & 2083 \\ \hline
%$\frac{1}{32}$& 55142 &32083\\ \hline
%\end{tabular}
%\caption{CPU time for Voronoi mesh.} % title name of the table
%\label{tab:e2h1err1-11t}
%\end{table}

%% file: ex2comp1.tex
\begin{table}[H]
\small
\centering % centering table
\begin{tabular}{| c| c c||  c c|} % creating 10 columns
\hline % inserting double-line
\multicolumn{5}{|c|}{ Iteration method.}\\
\hline
& \multicolumn{2}{|c||}{$u_1$} & \multicolumn{2}{|c|}{$u_2$} \\\cline{2-5}
$h$&  $e_{h,0}$  & $e_{h,1}$  & $e_{h,0}$  & $e_{h,1}$   \\ \hline
$\frac{1}{8}$& $5.839672e^{-4}$  & $6.380904e^{-2}$  &  $4.849473e^{-4}$ & $5.060191e^{-2}$ \\ \hline
$\frac{1}{16}$&$7.722887e^{-5}$  & $1.598596e^{-2}$  &  $6.083321e^{-5}$ & $1.264378e^{-2}$ \\ \hline
$\frac{1}{32}$& $2.739111e^{-5}$  & $4.000714e^{-3}$  & $1.014168e^{-5}$ & $3.159645e^{-3}$ \\ \hline
&\multicolumn{2}{|c||}{$u_3$} & \multicolumn{2}{|c|}{$u_4$} \\\cline{2-5}
$h$  &  $e_{h,0}$  & $e_{h,1}$  & $e_{h,0}$  & $e_{h,1}$   \\ \hline
$\frac{1}{8}$& $1.982376e^{-4}$  & $2.062463e^{-2}$  &  $3.430639e^{-3}$ & $3.807783e^{-1}$ \\ \hline
$\frac{1}{16}$&$2.786719e^{-5}$  & $5.174224e^{-3}$  &  $4.306165e^{-4}$ & $9.573263e^{-2}$ \\ \hline
$\frac{1}{32}$& $1.308594e^{-5}$  & $1.296721e^{-3}$  & $5.438818e^{-5}$ & $2.398772e^{-2}$ \\ \hline
\end{tabular}
\caption{Error values $e_{h,0},\,e_{h,1}$ for nonconvex mesh.} % title name of the table
\label{tab:e2h1err1-2}
\end{table}

\begin{table}[H]
\small
\centering % centering table
\begin{tabular}{| c| c c||  c c|} % creating 10 columns
\hline % inserting double-line
\multicolumn{5}{|c|}{ Two-grid method.}\\
\hline
& \multicolumn{2}{|c||}{$u_1$} & \multicolumn{2}{|c|}{$u_2$} \\\cline{2-5}
$h$&  $e_{h,0}$  & $e_{h,1}$  & $e_{h,0}$  & $e_{h,1}$   \\ \hline
$\frac{1}{8}$& $5.835024e^{-4}$  & $6.380899e^{-2}$  &  $4.848240e^{-4}$ & $5.060195e^{-2}$ \\ \hline
$\frac{1}{16}$&$7.370203e^{-5}$  & $1.598590e^{-2}$  &  $6.132023e^{-5}$ & $1.264375e^{-2}$ \\ \hline
$\frac{1}{32}$& $2.687501e^{-5}$  & $4.000648e^{-3}$  & $1.103731e^{-5}$ & $3.159606e^{-3}$ \\ \hline
&\multicolumn{2}{|c||}{$u_3$} & \multicolumn{2}{|c|}{$u_4$} \\\cline{2-5}
$h$  &  $e_{h,0}$  & $e_{h,1}$  & $e_{h,0}$  & $e_{h,1}$   \\ \hline
$\frac{1}{8}$& $1.980416e^{-4}$  & $2.062461e^{-2}$  &  $3.430072e^{-3}$ & $3.807783e^{-1}$ \\ \hline
$\frac{1}{16}$&$2.619106e^{-5}$  & $5.174197e^{-3}$  &  $4.353126e^{-4}$ & $9.573265e^{-2}$ \\ \hline
$\frac{1}{32}$& $1.293289e^{-5}$  & $1.296694e^{-3}$  & $5.470843e^{-5}$ & $2.398773e^{-2}$ \\ \hline
\end{tabular}
\caption{Error values $e_{h,0},\,e_{h,1}$ for nonconvex mesh.} % title name of the table
\label{tab:e2h1err1-22}
\end{table}

\begin{table}[H]
\small
\centering % centering table
\begin{tabular}{| c| c| c|   c||c|c|} % creating 10 columns
\hline
& &\multicolumn{2}{|c||}{Voronoi mesh} & \multicolumn{2}{|c|}{Nonconvex mesh}\\ \cline{3-6}
$h$& $H$& Iteration method & Two-grid method& Iteration method & Two-grid method \\
& &Time (sec)&Time(sec)& Time (sec)&Time(sec)\\ \hline
$\frac{1}{8}$& $\frac{1}{4}$& 387  &306&441&355\\ \hline
$\frac{1}{16}$&$\frac{1}{4}$& 2842  & 1846&3249&2131 \\ \hline
$\frac{1}{32}$&$\frac{1}{8}$&  55142 &29915&57954&30776\\ \hline
\end{tabular}
\caption{CPU time for Voronoi and nonconvex mesh.} % title name of the table
\label{tab:e2h1err1-11t}
\end{table}